\documentclass[12pt,letterpaper]{article}%
\usepackage{amsmath}
\usepackage{amsfonts}
\usepackage{amssymb}
\usepackage{graphicx}
\usepackage{color}%
\providecommand{\U}[1]{\protect\rule{.1in}{.1in}}
\newtheorem{theorem}{Theorem}[section]

\newtheorem{corollary}[theorem]{Corollary}

\newtheorem{definition}{Definition}

\newtheorem{lemma}[theorem]{Lemma}

\newtheorem{observation}[theorem]{Observation}
\newtheorem{problem}[theorem]{Problem}
\newtheorem{question}{Question}
\newtheorem{proposition}[theorem]{Proposition}
\newtheorem{remark}[theorem]{Remark}

\begin{document}

\title{\textbf{On }$k$-\textbf{Total Dominating Graphs}}
\author{S. Alikhani and D. Fatehi\\Department of Mathematics\\Yazd University, 89195-741, Yazd, \textsc{Iran}\\{\small alikhani@yazd.ac.ir, davidfatehi@yahoo.com}
\and C. M. Mynhardt\thanks{Supported by the Natural Sciences and Engineering
Research Council of Canada.}\\Department of Mathematics and Statistics\\University of Victoria, Victoria, BC, \textsc{Canada}\\{\small kieka@uvic.ca}}
\maketitle

\begin{abstract}
For a graph $G$, the $k$-total dominating graph $D_{k}^{t}(G)$ is the graph
whose vertices correspond to the total dominating sets of $G$ that have
cardinality at most $k$; two vertices of $D_{k}^{t}(G)$ are adjacent if and
only if the corresponding total dominating sets of $G$ differ by either adding
or deleting a single vertex. The graph $D_{k}^{t}(G)$ is used to study the
reconfiguration problem for total dominating sets: a total dominating set can
be reconfigured to another by a sequence of single vertex additions and
deletions, such that the intermediate sets of vertices at each step are total
dominating sets, if and only if they are in the same component of $D_{k}%
^{t}(G)$. Let $d_{0}(G)$ be the smallest integer $\ell$ such that $D_{k}%
^{t}(G)$ is connected for all $k\geq\ell$.

We investigate the realizability of graphs as total dominating graphs. For $k$
the upper total domination number $\Gamma_{t}(G)$, we show that any graph
without isolated vertices is an induced subgraph of a graph $G$ such that
$D_{k}^{t}(G)$ is connected. We obtain the bounds $\Gamma_{t}(G)\leq
d_{0}(G)\leq n$ for any connected graph $G$ of order $n\geq3$, characterize
the graphs for which either bound is realized, and determine $d_{0}(C_{n})$
and $d_{0}(P_{n})$.

\end{abstract}

\noindent\textbf{Keywords:} Total domination; Total domination reconfiguration
problem; $k$-Total dominating graph

\noindent\textbf{AMS Subject Classification Number 2010:\hspace{0.1in}}05C69

\section{Introduction}

A \emph{total dominating set} (\emph{TDS}) of a graph $G=(V,E)$ without
isolated vertices is a set $S\subseteq V(G)$ such that every vertex of $G$ is
adjacent to a vertex in $S$. If no proper subset of $S$ is a TDS of $G$, then
$S$ is \emph{a minimal TDS} (\emph{an MTDS}) of $G$. Every graph without
isolated vertices has a TDS, since $S=V(G)$ is such a set. The \emph{total
domination number} of $G$, denoted by $\gamma_{t}(G)$, is the minimum
cardinality of a TDS. The \emph{upper total domination number} of $G$, denoted
by $\Gamma_{t}(G)$, is the maximum cardinality of an MTDS. A TDS of size
$\gamma_{t}$ is called $\gamma_{t}$-\emph{set} of $G$, and an MTDS of size
$\Gamma_{t}(G)$ is called a $\Gamma_{t}$-\emph{set}.

For a given threshold $k$, let $S$ and $S^{\prime}$ be total dominating sets
of order at most $k$ of $G$. The \emph{total dominating set reconfiguration
(TDSR) problem} asks whether there exists a sequence of total dominating sets
of $G$ starting with $S$ and ending with $S^{\prime}$, such that each total
dominating set in the sequence is of order at most $k$ and can be obtained
from the previous one by either adding or deleting exactly one vertex. This
problem is similar to the \emph{dominating set reconfiguration (DSR) problem},
which is PSPACE-complete even for planar graphs, bounded bandwidth graphs,
split graphs, and bipartite graphs, while it can be solved in linear time for
cographs, trees, and interval graphs~\cite{haddadan-2015}.

The DSR problem naturally leads to the concept of the $k$-dominating graph
introduced by Haas and Seyffarth~\cite{HS1} as follows. If $G$ is a graph and
$k$ a positive integer, then the \emph{$k$-dominating graph} $D_{k}(G)$ of $G$
is the graph whose vertices correspond to the dominating sets of $G$ that have
cardinality at most $k$, two vertices of $D_{k}(G)$ being adjacent if and only
if the corresponding dominating sets of $G$ differ by either adding or
deleting a single vertex. The DSR problem therefore simply asks whether two
given vertices of $D_{k}(G)$ belong to the same component of $D_{k}(G)$. The
Haas-Seyffarth paper~\cite{HS1} stimulated the work of Alikhani, Fatehi and
Klav\v{z}ar \cite{davood}, Mynhardt, Roux and Teshima \cite{MRT}, Suzuki,
Mouawad and Nishimura \cite{SMN}, as well as their own follow-up paper
\cite{HS2}.

The study of $k$-dominating graphs was further motivated by similar studies of
graph colourings and other graph problems, such as independent sets, cliques
and vertex covers -- see e.g.~\cite{BC, CHJ1, CHJ2, Ito2, IKD} -- and by a
general goal to further understand the relationship between the dominating
sets of a graph. Motivated by definition of $k$-dominating graph, we define
the $k$-total dominating graph of $G$ as follows.

\begin{definition}
The \emph{$k$-total dominating graph} $D_{k}^{t}(G)$ of $G$ is the graph whose
vertices correspond to the total dominating sets of $G$ that have cardinality
at most $k$. Two vertices of $D_{k}^{t}(G)$ are adjacent if and only if the
corresponding total dominating sets of $G$ differ by either adding or deleting
a single vertex. For $r\geq0$, we abbreviate $D_{\Gamma_{t}(G)+r}^{t}(G)$ to
$D_{\Gamma_{t}+r}^{t}(G)$, and $D_{\gamma_{t}(G)+r}^{t}(G)$ to $D_{\gamma
_{t}+r}^{t}(G)$.
\end{definition}

In studying the TDSR problem, it is therefore natural to determine conditions
for $D_{k}^{t}(G)$ to be connected. We begin the study of this problem in
Section \ref{Sec_Connected}. To this purpose we define $d_{0}(G)$ to be the
smallest integer $\ell$ such that $D_{k}^{t}(G)$ is connected for all
$k\geq\ell$, and note that $d_{0}(G)$ exists for all graphs $G$ without
isolated vertices because $D_{|V(G)|}^{t}(G)$ is connected.

We introduce our notation in Section \ref{Sec_Not} and provide background
material on total domination in Section \ref{Sec_Back}. For instance, we
characterize graphs $G$ such that $\Gamma_{t}(G)=|V(G)|-1$. In Section
\ref{Sec_Cycle} we determine $d_{0}(C_{n})$ and $d_{0}(P_{n})$; interestingly,
it turns out that $d_{0}(C_{8})=\Gamma_{t}(C_{8})+2$, making $C_{8}$ the only
known graph for which $D_{\Gamma_{t}+1}^{t}(G)\ $is disconnected. The main
result for cycles requires four lemmas, which we state in Section
\ref{Sec_Cycle} but only prove in Section \ref{Sec_Proofs} to improve the flow
of the exposition. In Section \ref{Sec_Real} we study the realizability of
graphs as total dominating graphs. We show that the hypercubes $Q_{n}$ and
stars $K_{1,n}$ are realizable for all $n\geq2$, that $C_{4},C_{6},C_{8}$ and
$C_{10}$ are the only realizable cycles, and that $P_{1}$ and $P_{3}$ are the
only realizable paths. Section \ref{Sec_Pr} contains a list of open problems
and questions for future consideration.

\subsection{Notation}

\label{Sec_Not}For domination related concepts not defined here we refer the
reader to \cite{HHS}. The monograph \cite{HY} by Henning and Yeo is a valuable
resource on total domination.

For vertices $u,v$ of a graph $G$, we write $u\sim v$ if $u$ and $v$ are
adjacent. A vertex $v$ such that $u\sim v$ for all $u\in V(G)-\{v\}$ is a
\emph{universal vertex}. We refer to a vertex of $G$ of degree $1$ as a
\emph{leaf} and to the unique neighbour of a leaf in $G\ncong K_{2}$ as a
\emph{stem}, and denote the number of leaves and stems of $G$ by $\lambda(G)$
and $\sigma(G)$, respectively. As usual, for $u,v\in V(G)$, $d(u,v)$ denotes
the distance from $u$ to $v$.

A set of cardinality $n$ is also called an $n$-\emph{set}. A subset of
cardinality $k$ of a set $A$ is called a $k$\emph{-subset} of $A$. The
\emph{hypercube} $Q_{n}$ is the graph whose vertices are the $2^{n}$ subsets
of an $n$-set, where two vertices are adjacent if and only if one set is
obtained from the other by deleting a single element.

The disjoint union of $r$ copies of a graph $H$ is denoted by $rH$. The
\emph{corona }$G\circ K_{1}$\emph{ of a graph }$G$ is the graph obtained by
joining each vertex of $G$ to a new leaf. A \emph{generalized corona} of $G$
is a graph obtained by joining each vertex of $G$ to one or more new leaves.
For a graph $G$ and a subset $U$ of $V(G)$, we denote by $G[U]$ the subgraph
of $G$ induced by $U$.

\begin{remark}
\label{RemStem}The set of stems of a graph $G$ is a subset of any TDS of $G$,
otherwise some leaf is not totally dominated. Hence $\gamma_{t}(G)\geq
\sigma(G)$.
\end{remark}

The \emph{open neighbourhood} of a vertex $v$ is $N(v)=\{u\in V(G):u\sim v\}$
and the \emph{closed neighbourhood} of $v$ is $N[v]=N(v)\cup\{v\}$.\emph{ }Let
$S\subseteq V(G)$. The \emph{closed neighbourhood} of $S$ is $N[S]=\bigcup
_{s\in S}N[s]$. The \emph{open private neighbourhood of a vertex }$s\in
S$\emph{ relative to }$S$, denoted $\operatorname{OPN}(s,S)$, consists of all
vertices in the open neighbourhood of $s$ that do not belong to the open
neighbourhood of any $s^{\prime}\in S-\{s\}$, that is, $\operatorname{OPN}%
(s,S)=N(s)-\bigcup_{s^{\prime}\in S-\{s\}}N(s^{\prime})$. A vertex in
$\operatorname{OPN}(s,S)$ may belong to $S$, in which case it is called an
\emph{internal private neighbour of }$s$ \emph{relative to }$S$, or it may
belong to $V(G)-S$, in which case it is called an \emph{external private
neighbour of }$s$ \emph{relative to }$S$. The set of internal (external,
respectively) private neighbours of $s$ relative to $S$ are denoted by
$\operatorname{IPN}(s,S)$ ($\operatorname{EPN}(s,S)$, respectively). Hence
$\operatorname{OPN}(s,S)=\operatorname{IPN}(s,S)\cup\operatorname{EPN}(s,S)$.
These sets play an important role in determining whether a TDS is an MTDS or not.

\subsection{Preliminary results}

\label{Sec_Back}Cockayne, Dawes and Hedetniemi \cite{CDH} characterize minimal
total dominating sets as follows.

\begin{proposition}
\label{Prop_MTDS}\emph{\cite{CDH}\hspace{0.1in}}A TDS $S$ of a graph $G$ is an
MTDS if and only if $\operatorname{OPN}(s,S)\neq\varnothing$ for every $s\in
S$.
\end{proposition}

We restate Proposition \ref{Prop_MTDS} in a more convenient form for our purposes.

\begin{corollary}
\label{Cor_Stem}\emph{\hspace{0.1in}}Let $S$ be a TDS of a graph $G$, $H$ the
subgraph of $G$ consisting of the components of $G[S]$ of order at least $3$,
and $X$ the set of stems of $H$. Then $S$ is an MTDS if and only if
$\operatorname{EPN}(s,S)\neq\varnothing$ for each $s\in V(H)-X$.
\end{corollary}

\noindent\textbf{Proof.\hspace{0.1in}}A vertex $s^{\prime}\in S$ belongs to
$\operatorname{IPN}(s,S)$ if and only if $s^{\prime}\sim s$ and $\deg
_{G[S]}(s^{\prime})=1$. Therefore, if $G[\{s,s^{\prime}\}]$ is a $K_{2}$
component of $G[S]$, then $s^{\prime}\in\operatorname{IPN}(s,S)$ and
$s\in\operatorname{IPN}(s^{\prime},S)$. Further, if $x$ is a stem of $H$, then
$x$ is adjacent to a vertex $x^{\prime}\in S$ such that $\deg_{G[S]}%
(x^{\prime})=1$, hence $x^{\prime}\in\operatorname{IPN}(x,S)$. By Proposition
\ref{Prop_MTDS}, therefore, $S$ is an MTDS if and only if $\operatorname{OPN}%
(s,S)\neq\varnothing$ for every $s\in V(H)-X$. But for any $s\in V(H)-X$, each
vertex in $N(s)\cap S$ is adjacent to at least two vertices in $S$, hence
$\operatorname{IPN}(s,S)=\varnothing$, which implies that $\operatorname{OPN}%
(s,S)\neq\varnothing$ if and only if $\operatorname{EPN}(s,S)\neq\varnothing
$.~$\blacksquare$

\bigskip

Cockayne et al. \cite{CDH} also established an upper bound on the total
domination number, while Favaron and Henning \cite{FH} established an upper
bound on the upper total domination number.

\begin{proposition}
\label{TD_Bounds}If $G$ is a connected graph of order $n\geq3$, then

\begin{enumerate}
\item[$(i)$] \emph{\cite{CDH}}\hspace{0.1in}$\gamma_{t}(G)\leq\frac{2n}{3}$, and

\item[$(ii)$] \emph{\cite{FH}}\hspace{0.1in}$\Gamma_{t}(G)\leq n-1$;
furthermore, if $G$ has minimum degree $\delta\geq2$, then $\Gamma_{t}(G)\leq
n-\delta+1$, and this bound is sharp.
\end{enumerate}
\end{proposition}

We now characterize graphs $G$ such that $\Gamma_{t}(G)=|V(G)|-1$.

\begin{proposition}
\label{Prop_n-1}A connected graph $G$ of order $n\geq3$ satisfies $\Gamma
_{t}(G)=n-1$ if and only if $n$ is odd and $G$ is obtained from $\frac{n-1}%
{2}K_{2}$ by joining a new vertex to at least one vertex of each $K_{2}$.
\end{proposition}

\noindent\textbf{Proof.\hspace{0.1in}}It is clear that $\Gamma_{t}(G)=n-1$ for
any such graph $G$. For the converse, assume $G$ is a connected graph of order
$n\geq3$ such that $\Gamma_{t}(G)=n-1$ and let $S$ be a $\Gamma_{t}$-set of
$G$. Suppose $H$ is a component of $G[S]$ of order at least $3$. If
$\delta(H)\geq2$, then $H$ has no stems. If $\delta(H)=1$, then $H$ has at
least as many leaves as stems, so that $H$ has at least two vertices that are
not stems. In either case Corollary \ref{Cor_Stem} implies that at least two
vertices in $S$ has nonempty external private neighbourhoods, which is
impossible since $|S|=n-1$. Therefore each component of $H$ is a $K_{2}$.
Since only one vertex of $G$ does not belong to $S$ and $G$ is connected, the
result follows.~$\blacksquare$

\section{Connectedness of $D^{t}_{k}(G)$}

\label{Sec_Connected}Haas and Seyffarth~\cite{HS1} showed that $D_{\Gamma
(G)}(G)$ is disconnected whenever $E(G)\neq\varnothing$. In contrast, we show
that any graph without isolated vertices is an induced subgraph of a graph $G$
such that $D_{\Gamma_{t}}^{t}(G)$ is connected. We obtain the bounds
$\Gamma_{t}(G)\leq d_{0}(G)\leq n$ for any connected graph $G$ of order
$n\geq3$, and characterize graphs that satisfy equality in either bound.

We begin with some definitions and basic results. For $k\geq\gamma_{t}(G)$ and
$A,B$ total dominating sets of $G$ of cardinality at most $k$, we write
$A\overset{k}{\leftrightarrow}B$, or simply $A\leftrightarrow B$ if $k$ is
clear from the context, if there is a path in $D_{k}^{t}(G)$ connecting $A$
and $B$. The binary relation $\leftrightarrow$ is clearly symmetric and
transitive. If $A\subseteq B$ and $b\in B-A$, then\thinspace$b$ is adjacent to
a vertex in $A$ because $A$ is a TDS. Hence $A\cup\{b\}$ is a TDS. Repeating
this argument shows that $A\leftrightarrow B$. More generally, if $C$ is also
a TDS of cardinality at most $k$ such that $A\cup B\subseteq C$, then
$A\overset{k}{\leftrightarrow}C\overset{k}{\leftrightarrow}B$. Repeating the
same argument if $B$ is a TDS and $A\subseteq B$ is an MTDS shows that if
$A\overset{k}{\leftrightarrow}A^{\prime}$ for all MTDS's $A,A^{\prime}$ of
cardinality at most $k$, then $D_{k}^{t}(G)$ is connected. We state these
facts explicitly for referencing.

\begin{observation}
\label{Obs_conn}Let $A,B,C$ be total dominating sets of a graph $G$ of
cardinality at most $k\geq\gamma_{t}(G)$.

\begin{enumerate}
\item[$(i)$] If $A\subseteq B$, then $A\overset{k}{\leftrightarrow}B.$

\item[$(ii)$] If $A\cup B\subseteq C$, then $A\overset{k}{\leftrightarrow
}C\overset{k}{\leftrightarrow}B$.

\item[$(iii)$] If $A\overset{k}{\leftrightarrow}A^{\prime}$ for all MTDS's
$A,A^{\prime}$ of cardinality at most $k$, then $D_{k}^{t}(G)$ is connected.
\end{enumerate}
\end{observation}

\begin{figure}[ptb]%
\centering
\includegraphics[
height=2.8055in,
width=4.3111in
]%
{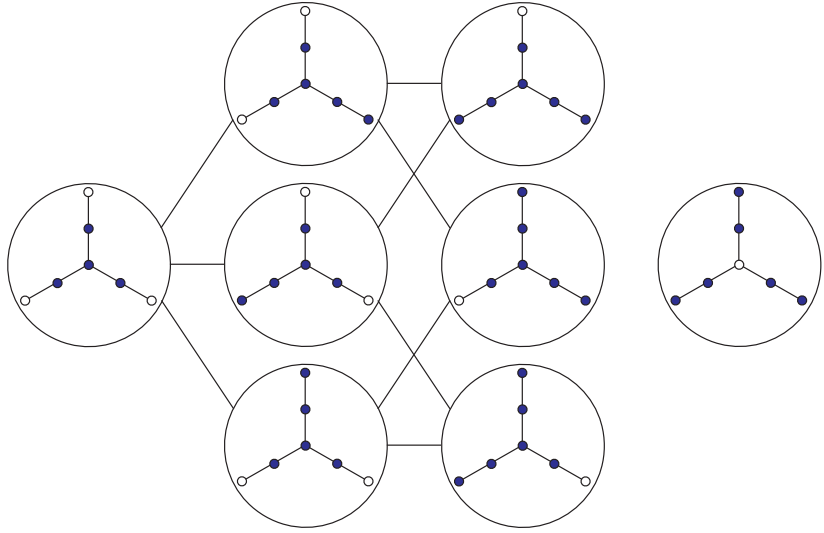}%
\caption{The graph $D_{6}^{t}(S_{2,2,2})$}%
\label{figure1}%
\end{figure}

As in the case of dominating sets, the connectedness of $D_{k}^{t}(G)$ however
does not guarantee the connectedness of $D_{k+1}^{t}(G)$. For example,
consider the tree $T=S_{2,2,2}$ (the spider with three legs of length $2$
each) in Figure \ref{figure1}. This figure shows $D_{6}^{t}(T)$, where
vertices are represented by copies of $T$, and the total dominating sets are
indicated by the solid circles. The unique $\Gamma_{t}$-set is an isolated
vertex in $D_{\Gamma}^{t}(T)$, so $D_{\Gamma_{t}}^{t}(T)=D_{n-1}^{t}(T)$ is
disconnected.%

In the case of dominating sets it is easy to see that $D_{\Gamma(G)}(G)$ is
disconnected whenever $G$ has at least one edge (and hence at least two
minimal dominating sets). For total domination the situation is not quite as
simple. A fundamental difference between domination and total domination is
that every graph with at least one edge has at least two different minimal
dominating sets, whereas there are many graphs with a unique MTDS. Consider,
for example, the double star $S(r,t)$, which consists of two adjacent vertices
$u$ and $v$ such that $u$ is adjacent to $r$ leaves and $v$ is adjacent to $t$
leaves. By Remark \ref{RemStem}, $u$ and $v$ belong to any TDS of $S(r,t)$.
Since $\{u,v\}$ is an MTDS, it is the \textbf{only} MTDS of $S(r,t)$.
Therefore $D_{\Gamma_{t}}^{t}(S(r,t))=K_{1}$, which is connected. We show that
the stems of $G$ determine whether $D_{\Gamma_{t}}^{t}(G)$ is connected or
not.\bigskip

\begin{theorem}
\label{ThmGammat_Conn}Let $G$ be a connected graph of order $n\geq3$ with
$\Gamma_{t}(G)=k$. Denote the set of stems of $G$ by $X$. Then $D_{k}^{t}(G)$
is connected if and only if $X$ is a TDS of $G$.
\end{theorem}

\noindent\textbf{Proof.\hspace{0.1in}}Let $S$ be any $\Gamma_{t}$-set of $G$.
Then no subset of $S$ is a TDS and no superset of $S$ is a vertex of
$D_{k}^{t}(G)$, hence $S$ is an isolated vertex of $D_{k}^{t}(G)$. Therefore
$D_{k}^{t}(G)$ is connected if and only if $S$ is the only MTDS of $G$.

Suppose $X$ is a TDS of $G$. Any $x\in X$ is adjacent to a leaf, hence
$X-\{x\}$ does not dominate $G$. Therefore $X$ is an MTDS. This implies that
no superset of $X$ is an MTDS. But by Remark \ref{RemStem}, $X$ is contained
in any TDS of $G$. Consequently, $X$ is the only MTDS of $G$, so $\gamma
_{t}(G)=\Gamma_{t}(G)=\sigma(G)$ and $D_{k}^{t}(G)=K_{1}$.

Conversely, suppose $X$ is not a TDS of $G$. We show that $G$ has at least two
MTDS's. First assume that $X$ dominates $G$. Then $G[X]$ has an isolated
vertex, say $x$, which is adjacent to a leaf $\ell_{x}\notin X$. Now%
\[
Y=X\cup\{\ell_{x}:x\ \text{is\ an\ isolated\ vertex\ of}\ G[X]\}
\]
is an MTDS of $G$. For another MTDS of $G$, let $T^{\prime}$ be a spanning
tree of $G$, let $T$ be the subtree of $T^{\prime}$ obtained by deleting all
leaves of $T^{\prime}$ and let $Z=V(T)$. If $|Z|=1$, then $T=K_{1}$ and
$T^{\prime}$ is a star. Say $Z=\{z\}$. Then $z$ is a universal vertex of $G$
and $X=\{z\}$. Since $n\geq3$ there exists a vertex $y\in N(z)-\ell_{z}$, and
$\{y,z\}$ is an MTDS of $G$ distinct from $Y$. On the other hand, if
$|Z|\geq2$, then $Z$ is a TDS that does not contain any leaves of $G$. Hence
$Z$ contains an MTDS distinct from $Y$.

Now assume that $X$ does not dominate $G$ and let $S$ be any MTDS of $G$. Then
there exists a vertex $v\in S-X$. Let $H$ be the component of $G[S]$ that
contains $v$ and consider the sets $\operatorname{IPN}(v,S)$ and
$\operatorname{EPN}(v,S)$. By Proposition \ref{Prop_MTDS}, $\operatorname{OPN}%
(v,S)\neq\varnothing$, hence at least one of $\operatorname{IPN}(v,S)$ and
$\operatorname{EPN}(v,S)$ is nonempty. Since $v\notin X$, $\deg_{G}(u)\geq2$
for each $u\in N(v)$. For each $u\in\operatorname{EPN}(v,S)$ we can therefore
choose a vertex $w_{u}\in V(G)-S$ adjacent to $u$, where possibly
$w_{u}=w_{u^{\prime}}$ for distinct $u,u^{\prime}\in\operatorname{EPN}(v,S)$.
Let $A=\operatorname{EPN}(v,S)\cup\{w_{u}:u\in\operatorname{EPN}(v,S)\}$.
Noting that if $u\in\operatorname{IPN}(v,S)$, then no other neighbour of $u$
belongs to $S$, we define the set $B=\operatorname{IPN}(v,S)\cup\{w_{u}%
:u\in\operatorname{IPN}(v,S)\}$ similarly. Let $S^{\prime}=(S-\{v\})\cup A\cup
B$. By definition $G[A]$ and $G[B]$ (if they are defined) have no isolated
vertices. If $B=\varnothing$, then $\operatorname{IPN}(v,S)=\varnothing$, and
as shown in the proof of Corollary \ref{Cor_Stem}, $H-v$ has no isolated
vertices. Therefore $G[S^{\prime}]$ has no isolated vertices. Since $A$ and
$B$ dominate $\operatorname{EPN}(v,S)$ and $\operatorname{IPN}(v,S)$,
respectively, $S^{\prime}$ dominates $G$. This shows that $S^{\prime}$ is a
TDS of $G$. Since $v\notin S^{\prime}$, $S^{\prime}$ contains an MTDS distinct
from $S$, which is what we wanted to show.~$\blacksquare$

\bigskip

The class of graphs whose stems form a TDS includes (but is not limited to)
the generalized coronas of graphs without isolated vertices. Hence any graph
without isolated vertices is an induced subgraph of a graph $G$ such that
$D_{\Gamma_{t}}^{t}(G)$ is connected. The first paragraph of the proof of
Theorem \ref{ThmGammat_Conn} implies the following result.

\begin{corollary}
\label{CorGammat_Conn}The graph $D_{\Gamma_{t}}^{t}(G)$ is disconnected if and
only if $G$ has at least two MTDS's.
\end{corollary}

If the set of stems of $G$ is a TDS, then it is the unique MTDS of $G$, hence
we also have the following corollary. The converse does not hold -- for the
spider $S(2,2,2)$ in Figure~\ref{figure1}, $D_{4}^{t}(S(2,2,2))=K_{1}$, which
is connected, but the stems form an independent set of cardinality $3$, which
is not a TDS.

\begin{corollary}
\label{CorGammat_Conn2}If $G$ is a connected graph of order $n\geq3$ whose set
of stems is a TDS, then $D_{\gamma_{t}}^{t}(G)$ is connected.
\end{corollary}

Since any TDS of cardinality greater than $\Gamma_{t}$ contains a TDS of
cardinality $\Gamma_{t}$, the following result is immediate from Observation
\ref{Obs_conn}$(i)$ (and similar to \cite[Lemma 4 ]{HS1}).

\begin{lemma}
\label{Lem_Con2}If $k\geq\Gamma_{t}(G)$ and $D_{k}^{t}(G)$ is connected, then
$D_{k+1}^{t}(G)$ is connected.
\end{lemma}

We now know that
\[
\Gamma_{t}(G)\leq d_{0}(G)\leq n
\]
for any connected graph of order $n\geq3$, and that the first inequality is
strict if and only if the stems of $G$ do not form a TDS. Equality in the
upper bound is realized by graphs with total domination number equal to $n-1$,
as characterized in Proposition \ref{Prop_n-1}, because all these graphs also
have an MTDS of cardinality $\frac{n-1}{2}+1$ different from the $\Gamma_{t}%
$-set described in the proof, so $D_{n-1}^{t}(G)$ is disconnected. We next
show that if $\Gamma_{t}(G)<n-1$, then $d_{0}(G)\leq n-1$.

\begin{theorem}
\label{Thm_less_n-1}If $G$ is a connected graph of order $n\geq3$ such that
$\Gamma_{t}(G)<n-1$, then $d_{0}(G)\leq\min\{n-1,\Gamma_{t}(G)+\gamma
_{t}(G)\}$.
\end{theorem}

\noindent\textbf{Proof.\hspace{0.1in}}Let $X$ be the set of stems of $G$.
Suppose first that $G$ has a unique MTDS $S$, so that $d_{0}(G)=\Gamma_{t}(G)$
by Corollary \ref{CorGammat_Conn}. By Remark \ref{RemStem}, $X$ is the unique
MTDS of $G$, hence $|X|\geq2$. But each vertex of $X$ is adjacent to a leaf,
hence $n\geq2|X|\geq|X|+2$. Therefore
\[
d_{0}(G)=\Gamma_{t}(G)\leq\left\lfloor \frac{n}{2}\right\rfloor \leq\left\{
\begin{tabular}
[c]{l}%
$n-1$\\
$\Gamma_{t}(G)+\gamma_{t}(G)$%
\end{tabular}
\ \right.  .
\]

Assume therefore that $G$ has at least two MTDS's and let $A$ and $B$ be any
two MTDS's of $G$. If $|A\cup B|\leq n-1$, then
$A\overset{n-1}{\leftrightarrow}B$ by Observation \ref{Obs_conn}$(ii)$, hence
assume $|A\cup B|=n$. By the hypothesis, $\Gamma_{t}(G)\leq n-2$, hence there
exist distinct vertices $a_{1},a_{2}\in A-B$ and $b_{1},b_{2}\in B-A$. By
Remark \ref{RemStem}, $\{a_{1},a_{2},b_{1},b_{2}\}\cap X=\varnothing$.
Consider the four pairs $a_{i},b_{j},\ i,j=1,2$. Suppose first that for one of
these pairs $a_{i},b_{j}$, every vertex adjacent to both $a_{i}$ and $b_{j}$
has degree at least $3$. Since we also have that $a_{i},b_{j}\notin X$,
$G-a_{i}-b_{j}$ has no isolated vertices. This implies that $V(G)-\{a_{i}%
,b_{j}\},\ V(G)-\{a_{i}\}$ and $V(G)-\{b_{j}\}$ are TDS's of $G$, and we have%
\[
A\overset{n-1}{\leftrightarrow}V(G)-\{b_{j}\}\overset{n-1}{\leftrightarrow
}V(G)-\{a_{i},b_{j}\}\overset{n-1}{\leftrightarrow}V(G)-\{a_{i}%
\}\overset{n-1}{\leftrightarrow}B.
\]
Hence assume that for each pair $a_{i},b_{j},\ i,j=1,2,$ there exists a vertex
$c_{i,j}$ such that $N(c_{i,j})=\{a_{i},b_{j}\}$. Then in $G_{1}%
=G-a_{1}-c_{1,1}$, $b_{1}\sim c_{2,1}$, and $c_{1,2}\sim b_{2}$. Since
$\deg(c_{1,1})=2$, $a_{1}$ and $c_{1,1}$ have no common neighbours except
possibly $b_{1}$ (which is adjacent to $c_{2,1}$). Therefore $G_{1}$ has no
isolated vertices, which means that $V(G)-\{a_{1},c_{1,1}\}$ is a TDS of $G$.
Similarly, $V(G)-\{b_{1},c_{1,1}\}$ and $V(G)-\{c_{1,1}\}$ are TDS's of $G$.
Now%
\[%
\begin{tabular}
[c]{ll}%
$A$ & $\overset{n-1}{\leftrightarrow}V(G)-\{b_{1}%
\}\overset{n-1}{\leftrightarrow}V(G)-\{b_{1},c_{1,1}%
\}\overset{n-1}{\leftrightarrow}V(G)-\{c_{1,1}\}$\\
& $\overset{n-1}{\leftrightarrow}V(G)-\{a_{1},c_{1,1}%
\}\overset{n-1}{\leftrightarrow}V(G)-\{a_{1}\}\overset{n-1}{\leftrightarrow
}B.$%
\end{tabular}
\
\]
By Observation \ref{Obs_conn}$(iii)$, $d_{0}(G)\leq n-1$.

Now let $C$ be any fixed $\gamma_{t}$-set and $B$ any MTDS of $G$. Then
$|C\cup B|\leq|C|+|B|\leq\gamma_{t}(G)+\Gamma_{t}(G)$. By Observation
\ref{Obs_conn}$(ii)$, $C\overset{\gamma_{t}+\Gamma_{t}}{\leftrightarrow}B$. By
transitivity, $A\overset{\gamma_{t}+\Gamma_{t}}{\leftrightarrow}B$ for all
MTDS's $A,B$ of $G$, so by Observation \ref{Obs_conn}$(iii)$, $d_{0}%
(G)\leq\Gamma_{t}(G)+\gamma_{t}(G)$.~$\blacksquare$

\bigskip

To summarise, in this section we showed that

\begin{itemize}
\item for any connected graph $G$ of order $n\geq3$,
\begin{equation}
\Gamma_{t}(G)\leq d_{0}(G)\leq n. \label{eq_d0}%
\end{equation}

\item The lower bound in (\ref{eq_d0}) is realized if and only if $G$ has
exactly one MTDS, i.e., if and only if the stems of $G$ form a TDS.

\item The upper bound in (\ref{eq_d0}) is realized if and only if $\Gamma
_{t}(G)=n-1$, i.e., if and only if $n$ is odd and $G$ is obtained from
$\frac{n-1}{2}K_{2}$ by joining a new vertex to at least one vertex of each
$K_{2}$.
\end{itemize}

\section{Determining $d_{0}$ for paths and cycles}

\label{Sec_Cycle}Our aim in this section is to show (in Theorem
\ref{Thm_Cycles}) that $d_{0}(C_{8})=\Gamma_{t}(C_{8})+2$ and $d_{0}%
(C_{n})=\Gamma_{t}(C_{n})+1$ if $n\neq8$. Similar techniques can be used to
show that $d_{0}(P_{2})=\Gamma_{t}(P_{2})=d_{0}(P_{4})=\Gamma_{t}(P_{4})=2$
and $d_{0}(P_{n})=\Gamma_{t}(P_{n})+1$ if $n=3$ or $n\geq5$. We need four
lemmas (Lemmas \ref{Lem_Cycle1} -- \ref{Lem_Cycle4}) to obtain the result for
cycles. To enhance the logical flow of the paper, we only state the lemmas in
this section and defer their proofs to Section \ref{Sec_Proofs}.\emph{ }

It is easy to determine the total domination numbers of paths and cycles.

\begin{observation}
\label{Ob_td_Paths}\emph{\cite[Observation 2.9]{HY}}\hspace{0.1in}For $n\geq
3$,%
\[
\gamma_{t}(P_{n})=\gamma_{t}(C_{n})=\left\{
\begin{tabular}
[c]{cl}%
$\frac{n}{2}+1$ & if $n\equiv2\ (\operatorname{mod}\ 4)$\\
& \\
$\left\lceil \frac{n}{2}\right\rceil $ & otherwise.
\end{tabular}
\right.
\]

\end{observation}

The upper total domination number for paths was determined by Dorbec, Henning
and McCoy \cite{Dorbec}.\ 

\begin{proposition}
\label{Prop_TD_Paths}\emph{\cite{Dorbec}}\hspace{0.1in}For any $n\geq2$,
$\Gamma_{t}(P_{n})=2\left\lfloor \frac{n+1}{3}\right\rfloor .$
\end{proposition}

The proof of the following proposition on the upper total domination number of
cycles can be found in the appendix.

\begin{proposition}
\label{Prop_Gt_Cn}For any $n\geq3$,
\[
\Gamma_{t}(C_{n})=\left\{
\begin{tabular}
[c]{ll}%
$2\left\lfloor \frac{n}{3}\right\rfloor $ & if $n\equiv2\ (\operatorname{mod}%
\ 6)$\\
& \\
$\left\lfloor \frac{2n}{3}\right\rfloor $ & otherwise.
\end{tabular}
\ \ \right.
\]

\end{proposition}

Let $C_{n}=(v_{0},v_{1},...,v_{n-1},v_{0})$. When discussing subsets of
$\{v_{0},v_{1},...,v_{n-1}\}$ the arithmetic in the subscripts is performed
modulo $n$. We mention some obvious properties of minimal total dominating
sets of $C_{n}$.

\begin{remark}
\label{Rem_Cycles}Let $S$ be an MTDS of $C_{n}$. Then

\begin{enumerate}
\item[$(i)$] each component of $C_{n}[S]$ is either $P_{2},\ P_{3}$ or $P_{4}$;

\item[$(ii)$] each $P_{3}$ or $P_{4}$ component is preceded and followed by
exactly two consecutive vertices of $C_{n}-S$;

\item[$(iii)$] $C_{n}-S$ does not contain three consecutive vertices of
$C_{n}$.
\end{enumerate}
\end{remark}

Using the next four lemmas, we show in Theorem \ref{Thm_Cycles} that, with the
single exception of $n=8$, $d_{0}(C_{n})=\Gamma_{t}(C_{n})+1$. We only state
the lemmas here; their proofs are given in Section \ref{Sec_Proofs}. The first
lemma concerns MTDS's that induce $P_{3}$ or $P_{4}$ components.

\begin{lemma}
\label{Lem_Cycle1}Let $n\geq10$.

\begin{enumerate}
\item[$(i)$] If $S$ is an MTDS such that $C_{n}[S]$ contains a $P_{4}$
component, then $S$ is connected, in $D_{|S|+1}^{t}(C_{n})$, to an MTDS
without $P_{4}$ components.

\item[$(ii)$] If $S$ is an MTDS such that $C_{n}[S]$ contains two consecutive
$P_{3}$ components, then $S$ is connected, in $D_{|S|+1}^{t}(C_{n})$, to an
MTDS with fewer $P_{3}$ components.

\item[$(iii)$] If $S$ is an MTDS such that $C_{n}[S]$ contains at least one
$P_{3}$ and at least one $P_{2}$ component but no $P_{4}$ components, then $S$
is connected, in $D_{\Gamma_{t}+1}^{t}(C_{n})$, to an MTDS that has no $P_{3}$ components.
\end{enumerate}
\end{lemma}

The next lemma concerns MTDS's that induce only $P_{2}$ components. For
brevity we refer to such an MTDS as a $P_{2}$-\emph{MTDS}. For a $P_{2}$-MTDS
$S$, each $P_{2}$ component is followed by one or two vertices not belonging
to $S$. We refer to these $P_{2}$ components as $P_{2}\overline{P}_{1}$ and
$P_{2}\overline{P}_{2}$ \emph{components}, respectively. An MTDS $S$ is called
a \emph{maximum }$P_{2}$-MTDS\emph{ }if $S$ is a $P_{2}$-MTDS of maximum cardinality.

\begin{lemma}
\label{Lem_Cycle2}Let $S$ be a $P_{2}$-MTDS of $C_{n},\ n\geq10$.

\begin{enumerate}
\item[$(i)$] $S$ is a maximum $P_{2}$-MTDS if and only if $C_{n}[S]$ has at
most two $P_{2}\overline{P}_{2}$ components.

\item[$(ii)$] If $C_{n}[S]$ has at least one $P_{2}\overline{P}_{1}$ component
and $S^{\prime}$ is any $P_{2}$-MTDS such that $|S|\leq|S^{\prime}|\leq|S|+2$,
then $S\overset{|S|+3}{\leftrightarrow}S^{\prime}$.
\end{enumerate}
\end{lemma}

We next consider $C_{n},\ n\equiv0\ (\operatorname{mod}\ 4)$.

\begin{lemma}
\label{Lem_Cycle3}Suppose $n\geq8$ and $n\equiv0\ (\operatorname{mod}\ 4)$;
say $n=4k$. (By Observation \ref{Ob_td_Paths}, $\gamma_{t}(C_{n})=2k$.) Then

\begin{enumerate}
\item[$(i)$] $D_{2k+1}^{t}(C_{n})$ is disconnected;

\item[$(ii)$] if $n\geq12$, then $C_{n}$ has a $P_{2}$-MTDS $S^{\ast}$ such
that $|S^{\ast}|=2k+2$ and $C_{n}[S^{\ast}]$ has four $P_{2}\overline{P}_{1}$ components;

\item[$(iii)$] all $\gamma_{t}$-sets belong to the same component of
$D_{2k+2}^{t}(C_{n})$ and all $P_{2}$-MTDS's of cardinality $2k$ or $2k+2$
belong to the same component of $D_{2k+3}^{t}(C_{n})$.
\end{enumerate}
\end{lemma}

Our final lemma concerns small cycles.

\begin{lemma}
\label{Lem_Cycle4}If $3\leq n\leq9$ and $n\neq8$, then $d_{0}(C_{n}%
)=\Gamma_{t}(C_{n})+1$.
\end{lemma}

\begin{theorem}
\label{Thm_Cycles}For $n=8$, $d_{0}(C_{8})=\Gamma_{t}(C_{8})+2$. In all other
cases, $d_{0}(C_{n})=\Gamma_{t}(C_{n})+1$.
\end{theorem}

\noindent\textbf{Proof.\hspace{0.1in}}Since $\gamma_{t}(C_{8})=\Gamma
_{t}(C_{8})=4$, Lemma \ref{Lem_Cycle3}$(i)$ implies that $D_{\Gamma_{t}+1}%
^{t}(C_{8})$ is disconnected and then the first part of Lemma \ref{Lem_Cycle3}%
$(iii)$ implies that $d_{0}(C_{8})=\Gamma_{t}(C_{8})+2$. By Lemma
\ref{Lem_Cycle4}, the theorem is true for $3\leq n\leq7$ and $n=9$. Hence
assume $n\geq10$.

Let $S$ be any MTDS of $C_{n}$. By Lemma \ref{Lem_Cycle1} (possibly applied
several times), if $C_{n}[S]$ has a $P_{3}$ or $P_{4}$ component, then there
exists a $P_{2}$-MTDS $S^{\ast}$ such that $S\overset{\Gamma_{t}%
+1}{\leftrightarrow}S^{\ast}$. Thus we may assume that $S$ is a $P_{2}$-MTDS.
If $n\not \equiv 0\ (\operatorname{mod}\ 4)$, then $S$ has at least one
$P_{2}\overline{P}_{1}$ component. If $n\equiv0\ (\operatorname{mod}\ 4)$,
then $n\geq12$ and, by Lemma \ref{Lem_Cycle3}, all $P_{2}$-MTDS's of
cardinality $\frac{n}{2}$ or $\frac{n}{2}+2$ belong to the same component of
$D_{2k+3}^{t}(C_{n})$. Moreover, any $P_{2}$-MTDS of cardinality $\frac{n}%
{2}+2$ has a $P_{2}\overline{P}_{1}$ component. In either case repeated
application of Lemma \ref{Lem_Cycle2}$(ii)$ shows that all $P_{2}$-MTDS's
belong to the same component of $D_{\Gamma_{t}+1}^{t}(C_{n})$. The result
follows from Observation \ref{Obs_conn}$(iii)$, Corollary \ref{CorGammat_Conn}
and Lemma \ref{Lem_Con2}.~$\blacksquare$

\bigskip

The proof of the following result for paths is similar and omitted. Note that
for $n\equiv0\ (\operatorname{mod}\ 4)$, $\Gamma_{t}(P_{n})=\Gamma_{t}%
(C_{n})+2$, which explains the difference between $d_{0}(P_{8})$ and
$d_{0}(C_{8})$. The result is trivial for $P_{2}=K_{2}$, while the result for
$P_{4}$ follows from Corollary \ref{CorGammat_Conn2}.

\begin{theorem}
\label{Thm_Paths}$d_{0}(P_{2})=\Gamma_{t}(P_{2})=d_{0}(P_{4})=\Gamma_{t}%
(P_{4})=2$ and $d_{0}(P_{n})=\Gamma_{t}(P_{n})+1$ if $n=3$ or $n\geq5$.
\end{theorem}

\section{Realizability of graphs as total dominating graphs}

\label{Sec_Real}One of the main problems in the study of $k$-total dominating
graphs is determining which graphs are total dominating graphs. Since
$D_{k}^{t}(G)=H$ if and only if $D_{k+2}^{t}(G\cup K_{2})=H$, in studying
graphs $G$ such that $D_{k}^{t}(G)=H$ for a given graph $H$ we restrict our
investigation to graphs $G$ without $K_{2}$ components (and also without
isolated vertices, so that $\gamma_{t}(G)$ is defined).

As noted in \cite{davood, MRT} for the $k$-dominating graph $D_{k}(G)$ of a
graph $G$ of order $n$, the $k$-total dominating graph $D_{k}^{t}(G)$ is
similarly a subgraph of the hypercube $Q_{n}$ (provided $k\geq\gamma_{t}(G)$
and $G$ has no isolated vertices) and is therefore bipartite. Since any subset
of $V(K_{n})$ of cardinality at least $2$ is a TDS of $K_{n}$ and since
$Q_{n}$ is vertex transitive, $D_{n}^{t}(K_{n})\cong Q_{n}-N[v]$ for some
$v\in V(Q_{n})$. We show in Corollary \ref{Cor_Qn_K1,n}$(i)$ that $Q_{n}$
itself is realizable as the $k$-total dominating graph of several graphs, and
in Corollary \ref{Cor_Qn_K1,n}$(ii)$ that stars $K_{1,n}$, $n\geq2$, are
realizable. Again the set of stems plays an important role.

In the last two results of the section we determine the realizability of paths
and cycles.

\begin{theorem}
\label{Thm_Qn}Let $H$ be any graph of order $r$, $2\leq r\leq n$, without
isolated vertices and let $G$ be a generalized corona of $H$ having exactly
$n$ leaves. For each $\ell$ such that $0\leq\ell\leq n$, $D_{r+\ell}^{t}(G)$
is the subgraph of $Q_{n}$ corresponding to the collection of all $k$-subsets,
$0\leq k\leq\ell$, of an $n$-set.
\end{theorem}

\noindent\textbf{Proof.\hspace{0.1in}}Every vertex of $H$ is a stem of $G$. By
Remark \ref{RemStem}, $X=V(H)$ is contained in any TDS of $G$. Since $H$ has
no isolated vertices, $X$ is an MTDS of $G$. As shown in the proof of Theorem
\ref{ThmGammat_Conn}, $X$ is the only MTDS of $G$. For any set $L$ of leaves
of $G$, $X\cup L$ is a TDS of $G$. Moreover, for any sets $L_{1}$ and $L_{2}$
of leaves, $X\cup L_{1}\overset{r+\ell}{\leftrightarrow}X\cup L_{2}$ if and
only if each $|L_{i}|\leq\ell$ and $L_{1}$ is obtained from $L_{2}$ by adding
or deleting exactly one vertex. The result now follows from the definitions of
$Q_{n}$ and $D_{r+\ell}^{t}(G)$.~$\blacksquare$

\begin{corollary}
\label{Cor_Qn_K1,n}Let $H$ be any graph of order $r$, $2\leq r\leq n$, without
isolated vertices and let $G$ be a generalized corona of $H$ having exactly
$n$ leaves. For every integer $n\geq2$,

\begin{enumerate}
\item[$(i)$] $D_{r+n}^{t}(G)\cong Q_{n}$

\item[$(ii)$] $D_{r+1}^{t}(G)\cong K_{1,n}$.
\end{enumerate}
\end{corollary}

\noindent\textbf{Proof.\hspace{0.1in}}$(i)\hspace{0.1in}$By Theorem
\ref{Thm_Qn}, $D_{r+n}^{t}(G)$ is the subgraph of $Q_{n}$ corresponding to the
collection of all subsets of an $n$-set. Hence $D_{r+n}^{t}(G)\cong Q_{n}$.

$(ii)\hspace{0.1in}$By Theorem \ref{Thm_Qn}, $D_{r+1}^{t}(G)$ is the subgraph
of $Q_{n}$ corresponding to the empty set and all singleton subsets of an
$n$-set. Hence $D_{r+1}^{t}(G)\cong K_{1,n}$.~$\blacksquare$

\bigskip

We mentioned above that for a graph $G$ without isolated vertices and
$\gamma_{t}(G)\leq k$, $D_{k}^{t}(G)$ is a subgraph of $Q_{n}$. The strategy
used in the proof of Theorem \ref{Thm_Qn} enables us to be a little more
specific in many cases.

\begin{proposition}
\label{Prop_Delta}Let $G$ be a connected graph of order $n\geq3$ having
$\sigma(G)$ stems. For any $k\geq\gamma_{t}(G)$, $D_{k}^{t}(G)$ is a subgraph
of $Q_{n-\sigma(G)}$.

In particular, $D_{n}^{t}(G)$ is a subgraph of $Q_{n-\sigma(G)}$ in which
$V(G)$ has degree $\Delta(D_{n}^{t}(G))=n-\sigma(G)=\Delta(Q_{n-\sigma(G)})$.
\end{proposition}

\noindent\textbf{Proof.\hspace{0.1in}}Let $X$ be the set of stems of $G$. By
Remark \ref{RemStem}, $X$ is contained in any TDS of $G$. Hence all TDS's of
$G$ are subsets of $V(G)$ that contain $X$, and there are $2^{n-\sigma(G)}$
such sets. This shows that $D_{k}^{t}(G)$ is a subgraph of $Q_{n-\sigma(G)}$
for any $k\geq\gamma_{t}(G)$.

Now consider $D_{n}^{t}(G)$. For $v\in V(G)$, $G-v$ has an isolated vertex if
and only if $v\in X$. Therefore $V(G)-\{u\}$ is a TDS of $G$ if and only if
$u\in V(G)-X$, which implies that $\deg(V(G))=n-\sigma(G)$ in $D_{n}^{t}(G)$.
Let $S$ be any TDS of $G$; necessarily, $X\subseteq S$. There are at most
$n-|S|$ supersets of $S$ of cardinality $|S|+1$ that are TDS's, and at most
$|S|-\sigma(G)$ subsets of $S$ of cardinality $|S|-1$ that are TDS's. Hence in
$D_{n}^{t}(G)$, $\deg(S)\leq n-|S|+|S|-\sigma(G)=n-\sigma(G)$.~$\blacksquare$

\bigskip

Concerning the realizability of cycles, it is easily seen that $D_{4}%
^{t}(P_{4})\cong C_{4}$, $D_{3}^{t}(C_{4})\cong D_{5}^{t}(P_{6})\cong C_{8}$,
$D_{4}^{t}(C_{5})\cong C_{10}$ and, if $G$ is the graph obtained by joining
two leaves of $K_{1,3}$, then $D_{3}^{t}(G)\cong D_{3}^{t}(K_{1,3})\cong
C_{6}$. We show that $C_{2r},\ r\in\{2,3,4,5\}$, are the only cycles
realizable as $k$-total domination graphs.

\begin{proposition}
\label{Prop_R_Cycle}$(i)\hspace{0.1in}$There is no graph $G$ of order $n>6$
such that $D_{k}^{t}(G)\cong C_{m}$ for some integer $k$.

$(ii)\hspace{0.1in}$For $m>10$, there is no graph $G$ such that $D_{k}%
^{t}(G)\cong C_{m}$ for some integer $k$.
\end{proposition}

\noindent\textbf{Proof.\hspace{0.1in}}$(i)\hspace{0.1in}$Suppose to the
contrary that $D_{k}^{t}(G)\cong C_{m}$. Let $S$ be a $\gamma_{t}$-set of $G$.
Then $\deg(S)=2$ in $D_{k}^{t}(G)$. Since each superset of $S$ is a TDS of
$G$, $S$ has exactly two supersets of cardinality $|S|+1$. This implies that
$n-|S|=2$, i.e., $n-\gamma_{t}(G)=2$. But we know that $\gamma_{t}(G)\leq
\frac{2n}{3}$ (Proposition \ref{TD_Bounds}$(i)$) and so $n\leq6$, which is a contradiction.

$(ii)\hspace{0.1in}$Now suppose that $D_{k}^{t}(G)\cong C_{m}$, where $m>10$.
By $(i)$, $G$ has order $n\leq6$ and $n-\gamma_{t}(G)=2$. Say $D_{k}^{t}(G)$
is the cycle $(S_{1},S_{2},...,S_{m},S_{1})$. Since $n=\gamma_{t}(G)+2$, each
$S_{i}$ has cardinality $\gamma_{t}(G),\ \gamma_{t}(G)+1$ or $\gamma_{t}(G)+2$.

First assume $|S_{i}|=\gamma_{t}(G)+2=n$ for some $i$. Then $n\geq4$ and we
also have $k=n$. Since $S_{i}$ has degree $2$ in $D_{k}^{t}(G)$, $G$ has
exactly two TDS's of cardinality $n-1$. By Remark \ref{RemStem}, $G$ has
exactly two vertices that are not stems. Since $n\geq4$ (and $G$ has no
$K_{2}$ components), $G$ consists of two stems and two leaves, i.e., $G=P_{4}%
$. But $D_{4}^{t}(P_{4})\cong C_{4}$, contradicting $m>10$.

We may therefore assume that $k=n-1$ and $n-2\leq|S_{i}|\leq n-1$ for each
$i$. But then, by definition of adjacency in $D_{k}^{t}(G)$, $|S_{i}|=n-1$ for
$\frac{m}{2}$ values of $i$ and $|S_{i}|=n-2$ for the other $\frac{m}{2}$
values of $i$. Since $m>10$, $V(G)$ has at least six subsets of cardinality
$n-1$, which implies that $n\geq6$. Therefore $n=6$, $\gamma_{t}(G)=4,\ k=5$
and $m=12$, and each of the six $5$-subsets of $V(G)$ is a TDS. By Remark
\ref{RemStem}, $G$ has no stems and hence no leaves. Let $v$ be a vertex of
$G$ such that $\deg(v)=\Delta(G)$. If $\deg(v)=5$, then $\{u,v\}$ is a TDS for
any $u\in V(G)-\{v\}$, which contradicts $\gamma_{t}(G)=4$. If $\deg(v)=4$,
let $u$ be the unique vertex nonadjacent to $v$ and let $w$ be any vertex
adjacent to $u$. Then $G[\{u,v,w\}]\cong P_{3}$, so $\{u,v,w\}$ is a TDS, also
a contradiction. If $\deg(v)=3$, let $u_{1}$ and $u_{2}$ be the vertices
nonadjacent to $v$. Since $G$ has no leaves, each $u_{i}$ is adjacent to a
vertex $w_{i}\in N(v)$. Hence $\{v,w_{1},w_{2}\}$ is a TDS of cardinality at
most $3$, again a contradiction. Therefore $G$ is $2$-regular. But if
$G=2K_{3}$, then $G$ has $\binom{3}{2}^{2}=9>6$ TDS's of cardinality $4$, and
if $G=C_{6}$, then any vertex of $D_{5}^{t}(G)$ corresponding to five
consecutive vertices of $C_{6}$ has degree $3$. With this final contradiction
the proof is complete.~$\blacksquare$

\bigskip

The realizability of paths is somewhat similar to that of cycles in that only
a small number of paths are $k$-total dominating graphs. Since $D_{2}%
^{t}(K_{2})\cong P_{1}$ and $D_{3}^{t}(P_{3})\cong D_{3}^{t}(P_{4})\cong
D_{4}^{t}(P_{5})\cong P_{3}$, $P_{1}$ and $P_{3}$ are realizable. Indeed, they
are the only realizable paths, as we show next.

\begin{proposition}
For $m\neq1,3$, there is no graph $G$ such that $D_{k}^{t}(G)\cong P_{m}$ for
some integer $k$.
\end{proposition}

\noindent\textbf{Proof.\hspace{0.1in}}Suppose $G$ is a graph of order $n$ such
that $D_{k}^{t}(G)\cong P_{m}$ for some integer $k$. Say $P_{m}=(S_{1}%
,S_{2},...,S_{m})$, where each $S_{i}$ is a TDS of $G$. It is easy to examine
all graphs of order at most $3$, hence assume $n\geq4$.

If $S_{1}$ is a $\gamma_{t}$-set of $G$, then exactly one superset of $S_{1}$
of cardinality $|S_{1}|+1$, namely $S_{2}$, is a TDS. Since every superset of
$S_{1}$ is a TDS, $\gamma_{t}(G)=n-1\leq\frac{2n}{3}$ (by Proposition
\ref{TD_Bounds}$(i)$), hence $n\leq3$, contrary to the assumption above. Thus
we may assume that $S_{1}$ and (similarly) $S_{m}$ have cardinality at least
$\gamma_{t}(G)+1$. Therefore $S_{i}$, for some $1<i<m$, is a $\gamma_{t}$-set.
Exactly as in the proof of Proposition \ref{Prop_R_Cycle} we obtain that
$n\leq6$, $\gamma_{t}=n-2$ and each $S_{i}$ has cardinality $\gamma
_{t}(G),\ \gamma_{t}(G)+1$ or $\gamma_{t}(G)+2$. If $|S_{i}|=\gamma_{t}(G)+2$,
then, as shown in the proof of Proposition \ref{Prop_R_Cycle}, $G=P_{4}$ and
$D_{4}^{t}(P_{4})\cong C_{4}\ncong P_{m}$. Therefore $|S_{1}|=|S_{m}%
|=\gamma_{t}(G)+1=n-1=k$. Now, $S_{2}$ is a $\gamma_{t}$-set of cardinality
$n-2$, hence $S_{1}$ and $S_{3}$ are the only supersets of $S_{2}$ of
cardinality $n-1$, and $S_{2}$, in turn, is the only subset of $S_{1}$ that is
a TDS. Therefore
\begin{equation}
\operatorname{OPN}(v,S_{1})=\varnothing\ \text{for\ exactly\ one\ vertex}%
\ v\in S_{1}. \label{eq_S1}%
\end{equation}
Let $G_{1}=G[S_{1}]$. Since $4\leq n\leq6$, $3\leq|S_{1}|\leq5$. Suppose
$G_{1}$ contains a triangle, say $(a_{1},a_{2},a_{3},a_{1})$. Then by
(\ref{eq_S1}) there exists a $b_{i}\in\operatorname{OPN}(a_{i},S_{1})$ for
$i=1,2$ (without loss of generality), where $b_{1}\neq b_{2}$ and
$\{b_{1},b_{2}\}\cap\{a_{1},a_{2},a_{3}\}=\varnothing$. Since $|S_{1}|=n-1$,
$b_{1}$ or $b_{2}$ belongs to $S_{1}$.

Say $b_{1}\in S_{1}$. Then $b_{1}\in\operatorname{IPN}(a_{1},S_{1})$ and we
also have from (\ref{eq_S1}) that $\operatorname{OPN}(b_{1},S_{1}%
)\neq\varnothing$ or $\operatorname{OPN}(a_{3},S_{1})\neq\varnothing$. But if
$\operatorname{OPN}(a_{3},S_{1})\neq\varnothing$, then $n=6$ and
$\{a_{1},a_{2},a_{3}\}$ is a TDS of $G$ of cardinality $n-3$, which is not the
case. Therefore there exists $c_{1}\in\operatorname{OPN}(b_{1},S_{1})$. But
then $V(G)=\{a_{1},a_{2},a_{3},b_{1},b_{2},c_{1}\}$ and $\{a_{1},a_{2}%
,b_{1}\}$ is a TDS of $G$, again a contradiction. We conclude that $G_{1}$ is triangle-free.

If $G_{1}$ contains an $r$-cycle, $r\geq4$, then (\ref{eq_S1}) implies that at
least $r-1$ vertices of the cycle have private neighbours not on the cycle.
But then $n\geq7$, a contradiction. Therefore $G_{1}$ is acyclic. If $G_{1}$
has two $K_{2}$ components, then, by the restrictions on the order of $G_{1}$,
$G_{1}=2K_{2}$ and $\operatorname{IPN}(v,S_{1})\neq\varnothing$ for each $v\in
S_{1}$, contrary to (\ref{eq_S1}).

Suppose $G_{1}$ has a path component $P_{r}=(u_{1},u_{2},...,u_{r})$, $r\geq
3$. Then neither leaf of $P_{r}$ has an internal private neighbour, so, by
(\ref{eq_S1}), one of them has an external private neighbour. Say $u_{1}$ has
external private neighbour $w_{1}$. Since $|S_{1}|=n-1$, $V(G)=S_{1}%
\cup\{w_{1}\}$, and since $w_{1}\in\operatorname{EPN}(u_{1},S_{1})$, $\deg
_{G}(w_{1})=1$. Now, if $S_{1}$ is disconnected, then $G_{1}=K_{2}\cup P_{r}$,
and since $\deg_{G}(w_{1})=1$, $K_{2}$ is also a component of $G$. Hence
$G=K_{2}\cup P_{4}$, so $D_{5}^{t}(G)=D_{3}^{t}(P_{4})=P_{3}$. On the other
hand, if $S_{1}$ is connected, then $G$ is isomorphic to $P_{4},\ P_{5}$ or
$P_{6}$, in which case $D_{n-1}^{t}(G)$ is $P_{3}$ or $C_{8}$.

Finally, suppose $G_{1}$ is a tree but not a path. Then $G_{1}$ has at least
three leaves. By (\ref{eq_S1}), two of them have external private neighbours,
contrary to $|S_{1}|=n-1$.~$\blacksquare$

\bigskip

A \emph{full subgraph }of $Q_{n}$ is a subgraph that corresponds to all
subsets of cardinality at least $k$ of an $n$-set, for some integer $k$ such
that $0\leq k\leq n$. In this section we showed that

\begin{itemize}
\item all full subgraphs of $Q_{n}$, $n\geq2$, are realizable as $k$-total
dominating graphs.

\item In particular, $Q_{n}$ and $K_{1,n}$ are realizable for all $n\geq2$.
\end{itemize}

We also showed that

\begin{itemize}
\item $C_{4},C_{6},C_{8}$ and $C_{10}$ are the only realizable cycles, and

\item $P_{1}$ and $P_{3}$ are the only realizable paths.
\end{itemize}

\section{Proofs of lemmas in Section \ref{Sec_Cycle}}

\label{Sec_Proofs}This section contains the proofs of the lemmas stated in
Section \ref{Sec_Cycle}. To simplify our discussion of total dominating sets
of $C_{n}$, we encode each TDS $S$ using an $n$-tuple (or part of an
$n$-tuple) of the symbols $\circ$ and $\bullet$, where $\bullet$ in position
$i$ indicates that $v_{i-1}\in S$, while $\circ$ in position $i$ indicates
that $v_{i-1}\notin S$. For example, the MTDS $S=\{v_{0},v_{1},v_{4},v_{5}\}$
of $C_{8}$ is written as $S=(\bullet\bullet\circ\circ\bullet\bullet\circ
\circ)$. By Remark \ref{Rem_Cycles}$(ii)$, every $P_{3}$ or $P_{4}$ component
belongs to a code of the form $(\cdots\circ\circ\bullet\bullet\bullet
\circ\circ\cdots)$ or $(\cdots\circ\circ\bullet\bullet\bullet\bullet\circ
\circ\cdots)$, respectively. The $n$-tuples are often compressed by writing
the number of consecutive occurrences of $\circ$ or $\bullet$ above the
symbol; for example, we may write $(\cdots\circ\circ\bullet\bullet\bullet
\circ\circ\cdots)$ and $(\cdots\circ\circ\bullet\bullet\bullet\bullet
\circ\circ\cdots)$ as $(\cdots\overset{2}{\circ}\overset{3}{\bullet
}\overset{2}{\circ}\cdots)$ and $(\cdots\overset{2}{\circ}\overset{4}{\bullet
}\overset{2}{\circ}\cdots)$, respectively. When a $P_{2}$-component of $S$ can
be followed by one of two vertices of $C_{n}-S$, we write $(\cdots
\overset{2}{\bullet}\circ\cdots)$ without indicating a number above $\circ$.

When adding vertices to a TDS of $C_{n}$, for example, when adding a vertex to
a TDS $(\cdots\overset{2}{\circ}\overset{2}{\bullet}\overset{2}{\circ}\cdots)$
to form a TDS $(\cdots\overset{2}{\circ}\overset{3}{\bullet}\overset{1}{\circ
}\cdots)$, we write $(\cdots\overset{2}{\circ}\overset{2}{\bullet
}\overset{2}{\circ}\cdots)\rightarrow(\cdots\overset{2}{\circ}%
\overset{3}{\blacklozenge}\overset{1}{\lozenge}\cdots)$ to emphasize the
position of the addition. Conversely, when deleting a vertex from a TDS
$(\cdots\overset{2}{\circ}\overset{3}{\bullet}\overset{1}{\circ}\cdots)$ (say)
to form a TDS $(\cdots\overset{2}{\circ}\overset{2}{\bullet}\overset{2}{\circ
}\cdots)$, we write $(\cdots\overset{2}{\circ}\overset{3}{\bullet
}\overset{1}{\circ}\cdots)\rightarrow(\cdots\overset{2}{\circ}%
\overset{2}{\blacklozenge}\overset{2}{\lozenge}\cdots)$.

We restate the lemmas for convenience.

\bigskip

\noindent\textbf{Lemma \ref{Lem_Cycle1}\hspace{0.1in}}\emph{Let }$n\geq
10$\emph{. }

\begin{enumerate}
\item[$(i)$] \emph{If }$S$\emph{ is an MTDS such that }$C_{n}[S]$\emph{
contains a }$P_{4}$\emph{ component, then }$S$\emph{ is connected, in
}$D_{|S|+1}^{t}(C_{n})$\emph{, to an MTDS without }$P_{4}$\emph{ components.}

\item[$(ii)$] \emph{If }$S$\emph{ is an MTDS such that }$C_{n}[S]$\emph{
contains two consecutive }$P_{3}$\emph{ components, then }$S$\emph{ is
connected, in }$D_{|S|+1}^{t}(C_{n})$\emph{, to an MTDS with fewer }$P_{3}%
$\emph{ components.}

\item[$(iii)$] \emph{If }$S$\emph{ is an MTDS such that }$C_{n}[S]$\emph{
contains at least one }$P_{3}$\emph{ and at least one }$P_{2}$\emph{ component
but no }$P_{4}$\emph{ components, then }$S$\emph{ is connected, in }%
$D_{\Gamma_{t}+1}^{t}(C_{n})$\emph{, to an MTDS that has no }$P_{3}$\emph{
components. }
\end{enumerate}

\noindent\textbf{Proof.\hspace{0.1in}}$(i)\hspace{0.1in}$The result is easy to
see for $C_{6}$, hence assume $n\geq8$. By Remark \ref{Rem_Cycles}$(i)$, $S$
is of the form $(\cdots\overset{2}{\circ}\overset{4}{\bullet}\overset{2}{\circ
}\cdots)$. Consider the TDS $S^{\prime}$ with $|S^{\prime}|=|S|$ obtained by%
\[
(\cdots\overset{2}{\circ}\overset{4}{\bullet}\overset{2}{\circ}\cdots
)\rightarrow(\cdots\overset{2}{\circ}\overset{5}{\blacklozenge}%
\overset{1}{\lozenge}\cdots)\rightarrow(\cdots\overset{2}{\circ}%
\overset{2}{\blacklozenge}\overset{1}{\lozenge}\overset{2}{\blacklozenge
}\overset{1}{\circ}\cdots)=S^{\prime}.
\]
If $S^{\prime}$ is an MTDS, let $S^{\prime\prime}=S^{\prime}$. If $S^{\prime}$
is not an MTDS, then $S^{\prime}$ is of the form $(\cdots\overset{2}{\circ
}\overset{2}{\bullet}\overset{1}{\circ}\overset{2}{\bullet}\overset{1}{\circ
}\overset{3}{\bullet}\cdots)$ or $(\cdots\overset{2}{\circ}\overset{2}{\bullet
}\overset{1}{\circ}\overset{2}{\bullet}\overset{1}{\circ}\overset{4}{\bullet
}\cdots)$. In the former case, let $S^{\prime\prime}=(\cdots\overset{2}{\circ
}\overset{2}{\bullet}\overset{1}{\circ}\overset{2}{\bullet}%
\overset{2}{\lozenge}\overset{2}{\blacklozenge}\cdots)$ and in the later case
let $S^{\prime\prime}=(\cdots\overset{2}{\circ}\overset{2}{\bullet
}\overset{1}{\circ}\overset{2}{\bullet}\overset{2}{\lozenge}%
\overset{3}{\blacklozenge}\cdots)$. In all cases, $S^{\prime\prime}$ is an
MTDS having fewer $P_{4}$ components than $S$ such that
$S\overset{|S|+1}{\leftrightarrow}S^{\prime\prime}$ and $|S^{\prime\prime
}|\leq|S|$. The result follows by repeating this procedure.

$(ii)\hspace{0.1in}$The result follows from the operations%
\[
(\cdots\overset{2}{\circ}\overset{3}{\bullet}\overset{2}{\circ}%
\overset{3}{\bullet}\overset{2}{\circ}\cdots)\rightarrow(\cdots
\overset{2}{\circ}\overset{3}{\bullet}\overset{1}{\lozenge}%
\overset{4}{\blacklozenge}\overset{2}{\circ}\cdots)\rightarrow(\cdots
\overset{2}{\circ}\overset{2}{\blacklozenge}\overset{2}{\lozenge
}\overset{4}{\bullet}\overset{2}{\circ}\cdots)\rightarrow(\cdots
\overset{2}{\circ}\overset{2}{\bullet}\overset{1}{\lozenge}%
\overset{5}{\blacklozenge}\overset{2}{\circ}\cdots)\rightarrow(\cdots
\overset{2}{\circ}\overset{2}{\bullet}\overset{1}{\circ}%
\overset{2}{\blacklozenge}\overset{1}{\lozenge}\overset{2}{\blacklozenge
}\overset{2}{\circ}\cdots).
\]

$(iii)\hspace{0.1in}$If $C_{n}[S]$ contains at least one $P_{3}$ and at least
one $P_{2}$ component but no $P_{4}$ components, then $C_{n}[S]$ contains
either (a) two consecutive $P_{3}$ components or (b) a $P_{3}$ component
preceded and followed by a $P_{2}$ component. In the former case, $C_{n}[S]$
is of the form $(\cdots\overset{2}{\circ}\overset{3}{\bullet}\overset{2}{\circ
}\overset{3}{\bullet}\overset{2}{\circ}\cdots)$ and, by $(ii)$,
$S\overset{|S|+1}{\leftrightarrow}S^{\prime}$ where $S^{\prime}=(\cdots
\overset{2}{\circ}\overset{2}{\bullet}\overset{1}{\circ}\overset{2}{\bullet
}\overset{1}{\circ}\overset{2}{\bullet}\overset{2}{\circ}\cdots)$. Moreover,
$|S^{\prime}|=|S|$. In the latter case, $C_{n}[S]$ is of the form
$(\cdots\circ\overset{2}{\bullet}\overset{2}{\circ}\overset{3}{\bullet
}\overset{2}{\circ}\overset{2}{\bullet}\circ\cdots)$, and $S^{\prime\prime
}=(\cdots\circ\overset{2}{\bullet}\overset{1}{\circ}\overset{2}{\bullet
}\overset{1}{\circ}\overset{2}{\bullet}\overset{1}{\circ}\overset{2}{\bullet
}\circ\cdots)$ is an MTDS of larger cardinality having fewer $P_{3}$
components than $S$. Hence $S$ is not a $\Gamma_{t}$-set. The operations%
\[
(\cdots\circ\overset{2}{\bullet}\overset{2}{\circ}\overset{3}{\bullet
}\overset{2}{\circ}\overset{2}{\bullet}\circ\cdots)\rightarrow(\cdots
\circ\overset{2}{\bullet}\overset{1}{\lozenge}\overset{5}{\blacklozenge
}\overset{1}{\lozenge}\overset{2}{\bullet}\circ\cdots)\rightarrow(\cdots
\circ\overset{2}{\bullet}\overset{1}{\circ}\overset{2}{\blacklozenge
}\overset{1}{\lozenge}\overset{2}{\blacklozenge}\overset{1}{\circ
}\overset{2}{\bullet}\circ\cdots)=S^{\prime\prime}%
\]
show that $S\overset{|S^{\prime\prime}|+1}{\leftrightarrow}S^{\prime}$, hence
$S\overset{\Gamma_{t}+1}{\leftrightarrow}S^{\prime\prime}$. By repeating the
operations for (a) and (b) as necessary we obtain the desired
result.~$\blacksquare$

\bigskip

\noindent\textbf{Lemma \ref{Lem_Cycle2}\hspace{0.1in}}\emph{Let }$S$\emph{ be
a }$P_{2}$\emph{-MTDS of }$C_{n},\ n\geq10$\emph{.}

\begin{enumerate}
\item[$(i)$] $S$\emph{ is a maximum }$P_{2}$\emph{-MTDS if and only if }%
$C_{n}[S]$\emph{ has at most two }$P_{2}\overline{P}_{2}$\emph{ components.}

\item[$(ii)$] \emph{If }$C_{n}[S]$\emph{ has at least one }$P_{2}\overline
{P}_{1}$\emph{ component and }$S^{\prime}$\emph{ is any }$P_{2}$\emph{-MTDS
such that }$|S|\leq|S^{\prime}|\leq|S|+2$\emph{, then }%
$S\overset{|S|+3}{\leftrightarrow}S^{\prime}$\emph{.}
\end{enumerate}

\noindent\textbf{Proof.\hspace{0.1in}}$(i)$ Suppose $C_{n}[S]$ has $p$
components, $q$ of which are $P_{2}\overline{P}_{2}$ components. Then $|S|=2p$
and $n=4q+3(p-q)$. The result follows by comparing these numbers to the
formula for $\Gamma_{t}(C_{n})$ in Proposition \ref{Prop_Gt_Cn}.

$(ii)\hspace{0.1in}$First note that if $S$ is of the form $(\cdots
\circ\overset{2}{\bullet}\overset{1}{\circ}\overset{2}{\bullet}%
\overset{2}{\circ}\cdots\overset{2}{\bullet}\overset{2}{\circ}%
\overset{2}{\bullet}\overset{1}{\circ}\cdots)$, then repeating the operations%
\[
(\cdots\circ\overset{2}{\bullet}\overset{1}{\circ}\overset{2}{\bullet
}\overset{2}{\circ}\cdots\overset{2}{\bullet}\overset{2}{\circ}%
\overset{2}{\bullet}\overset{1}{\circ}\cdots)\rightarrow(\cdots\circ
\overset{2}{\bullet}\overset{1}{\circ}\overset{3}{\blacklozenge}%
\overset{1}{\lozenge}\cdots\overset{2}{\bullet}\overset{2}{\circ
}\overset{2}{\bullet}\overset{1}{\circ}\cdots)\rightarrow(\cdots
\circ\overset{2}{\bullet}\overset{2}{\lozenge}\overset{2}{\blacklozenge
}\overset{1}{\circ}\cdots\overset{2}{\bullet}\overset{2}{\circ}%
\overset{2}{\bullet}\overset{1}{\circ}\cdots)
\]
as necessary shows that $S$ is connected in $D_{|S|+1}^{t}(C_{n})$ to an MTDS
of the same cardinality, hence with the same number of both types of
components, in which all the components of each type appear consecutively.

Thus, assume $S$ is of the form $(\underline{\bullet}\bullet\overset{1}{\circ
}\overset{2}{\bullet}\overset{1}{\circ}\cdots\overset{2}{\bullet
}\overset{1}{\circ}\overset{2}{\bullet}\overset{2}{\circ}\overset{2}{\bullet
}\overset{2}{\circ}\cdots)$, where all the components of each type appear
consecutively and $\underline{\bullet}$ is a marked vertex to indicate the
position of the first $P_{2}\overline{P}_{1}$ component. The operations
\begin{align*}
(\underline{\bullet}\bullet\overset{1}{\circ}\overset{2}{\bullet
}\overset{1}{\circ}\overset{2}{\bullet}\overset{1}{\circ}\cdots
\overset{2}{\bullet}\overset{1}{\circ}\overset{2}{\bullet}\overset{2}{\circ
}\cdots\overset{2}{\bullet}\overset{2}{\circ})  &  \rightarrow
(\underline{\bullet}\overset{4}{\blacklozenge}\overset{1}{\circ}%
\overset{2}{\bullet}\overset{1}{\circ}\cdots\overset{2}{\bullet}%
\overset{1}{\circ}\overset{2}{\bullet}\overset{2}{\circ}\cdots
\overset{2}{\bullet}\overset{2}{\circ})\rightarrow(\underline{\lozenge
}\overset{4}{\bullet}\overset{1}{\circ}\overset{2}{\bullet}\overset{1}{\circ
}\cdots\overset{2}{\bullet}\overset{1}{\circ}\overset{2}{\bullet
}\overset{2}{\circ}\cdots\overset{2}{\bullet}\overset{2}{\circ})\\
&  \rightarrow(\underline{\circ}\overset{7}{\blacklozenge}\overset{1}{\circ
}\cdots\overset{2}{\bullet}\overset{1}{\circ}\overset{2}{\bullet
}\overset{2}{\circ}\cdots\overset{2}{\bullet}\overset{2}{\circ})\rightarrow
(\underline{\circ}\overset{2}{\blacklozenge}\overset{1}{\lozenge
}\overset{4}{\blacklozenge}\overset{1}{\circ}\cdots\overset{2}{\bullet
}\overset{1}{\circ}\overset{2}{\bullet}\overset{2}{\circ}\cdots
\overset{2}{\bullet}\overset{2}{\circ})\\
&  \rightarrow\cdots\rightarrow(\underline{\circ}\overset{2}{\bullet
}\overset{1}{\circ}\overset{2}{\bullet}\overset{1}{\circ}\overset{2}{\bullet
}\cdots\overset{4}{\blacklozenge}\overset{2}{\circ}\cdots\overset{2}{\bullet
}\overset{2}{\circ})\\
&  \rightarrow(\underline{\circ}\overset{2}{\bullet}\overset{1}{\circ
}\overset{2}{\bullet}\overset{1}{\circ}\overset{2}{\bullet}\cdots
\overset{5}{\blacklozenge}\overset{1}{\lozenge}\cdots\overset{2}{\bullet
}\overset{2}{\circ})\rightarrow(\underline{\circ}\overset{2}{\bullet
}\overset{1}{\circ}\overset{2}{\bullet}\overset{1}{\circ}\overset{2}{\bullet
}\cdots\overset{2}{\blacklozenge}\overset{1}{\lozenge}%
\overset{2}{\blacklozenge}\overset{1}{\circ}\cdots\overset{2}{\bullet
}\overset{2}{\circ})\\
&  \rightarrow\cdots\rightarrow(\underline{\circ}\overset{2}{\bullet
}\overset{1}{\circ}\overset{2}{\bullet}\overset{1}{\circ}\overset{2}{\bullet
}\cdots\overset{2}{\bullet}\overset{1}{\circ}\overset{2}{\bullet
}\overset{1}{\circ}\cdots\overset{2}{\lozenge}\overset{2}{\blacklozenge
}\overset{1}{\lozenge})
\end{align*}
produce a $P_{2}$-MTDS $S^{\prime\prime}$ which can also be obtained from $S$
by a rotation $v_{i}\rightarrow v_{i+1}$ for each $i$. Thus
$S\overset{|S|+1}{\leftrightarrow}S^{\prime\prime}$. These operations can be
repeated to show that $S\overset{|S|+1}{\leftrightarrow}S_{3}$ for each
rotation $S_{3}$ of $S$. By the above and transitivity, for each $P_{2}$-MTDS
$S^{\prime}$ such that $|S|=|S^{\prime}|$, $S\overset{|S|+1}{\leftrightarrow
}S^{\prime}$ and hence $S\overset{|S|+3}{\leftrightarrow}S^{\prime}$.

Now assume that $S^{\prime}$ is any $P_{2}$-MTDS such that $|S^{\prime
}|=|S|+2$. Then $S$ is not a maximum $P_{2}$-MTDS and hence, by $(i)$, $S$ has
at least three $P_{2}\overline{P}_{2}$ components. As shown above we may
assume all $P_{2}\overline{P}_{2}$ components of $C_{n}[S]$ are consecutive.
Hence $S$ is of the form $(\overset{2}{\bullet}\overset{2}{\circ
}\overset{2}{\bullet}\overset{2}{\circ}\overset{2}{\bullet}\overset{2}{\circ
}\cdots\overset{2}{\bullet}\overset{1}{\circ}\cdots)$, and the addition of
three vertices in succession produces a TDS of the form $(\overset{2}{\bullet
}\overset{1}{\lozenge}\overset{7}{\blacklozenge}\overset{2}{\circ}%
\cdots\overset{2}{\bullet}\overset{1}{\circ}\cdots)$. Then the operations%
\[
(\overset{2}{\bullet}\overset{1}{\circ}\overset{7}{\bullet}\overset{2}{\circ
}\cdots\overset{2}{\bullet}\overset{1}{\circ}\cdots)\rightarrow
(\overset{2}{\bullet}\overset{1}{\circ}\overset{2}{\blacklozenge
}\overset{1}{\lozenge}\overset{4}{\blacklozenge}\overset{2}{\circ}%
\cdots\overset{2}{\bullet}\overset{1}{\circ}\cdots)\rightarrow
(\overset{2}{\bullet}\overset{1}{\circ}\overset{2}{\bullet}\overset{1}{\circ
}\overset{5}{\blacklozenge}\overset{1}{\lozenge}\cdots\overset{2}{\bullet
}\overset{1}{\circ}\cdots)\rightarrow(\overset{2}{\bullet}\overset{1}{\circ
}\overset{2}{\bullet}\overset{1}{\circ}\overset{2}{\blacklozenge
}\overset{1}{\lozenge}\overset{2}{\blacklozenge}\overset{1}{\lozenge}%
\cdots\overset{2}{\bullet}\overset{1}{\circ}\cdots)=S_{1}%
\]
produce a $P_{2}$-MTDS $S_{1}$ such that $|S_{1}|=|S|+2$ and $S_{1}%
\overset{|S|+3}{\leftrightarrow}S$. However, we have already shown above that
$S_{1}\overset{|S_{1}|+1}{\leftrightarrow}S^{\prime}$, i.e. $S_{1}%
\overset{|S|+3}{\leftrightarrow}S^{\prime}$. By transitivity,
$S\overset{|S|+3}{\leftrightarrow}S^{\prime}$.~$\blacksquare$

\bigskip

\noindent\textbf{Lemma \ref{Lem_Cycle3}\hspace{0.1in}}\emph{Suppose }$n\geq
8$\emph{ and }$n\equiv0\ (\operatorname{mod}\ 4)$\emph{; say }$n=4k$\emph{.
(By Observation \ref{Ob_td_Paths}, }$\gamma_{t}(C_{n})=2k$\emph{.) Then}

\begin{enumerate}
\item[$(i)$] $D_{2k+1}^{t}(C_{n})$\emph{ is disconnected;}

\item[$(ii)$] \emph{if }$n\geq12$\emph{, then }$C_{n}$\emph{ has a }$P_{2}%
$\emph{-MTDS }$S^{\ast}$\emph{ such that }$|S^{\ast}|=2k+2$\emph{ and }%
$C_{n}[S^{\ast}]$\emph{ has four }$P_{2}\overline{P}_{1}$\emph{ components;}

\item[$(iii)$] \emph{all }$\gamma_{t}$\emph{-sets belong to the same component
of }$D_{2k+2}^{t}(C_{n})$\emph{ and all }$P_{2}$\emph{-MTDS's of cardinality
}$2k$\emph{ or }$2k+2$\emph{ belong to the same component of }$D_{2k+3}%
^{t}(C_{n})$\emph{.}
\end{enumerate}

\noindent\textbf{Proof.\hspace{0.1in}}Any\textbf{ }$\gamma_{t}$-set $S$ of
$C_{n}$ is a $P_{2}$-MTDS, hence of the form $(\cdots\overset{2}{\circ
}\overset{2}{\bullet}\overset{2}{\circ}\overset{2}{\bullet}\overset{2}{\circ
}\cdots)$.

$(i)\hspace{0.1in}$By symmetry the addition of any single vertex $v$ to $S$
results in $S^{\prime}=(\cdots\overset{2}{\circ}\overset{3}{\blacklozenge
}\overset{1}{\lozenge}\overset{2}{\bullet}\overset{2}{\circ}\cdots)$, and by
Remark \ref{Rem_Cycles}$(iii)$, $v$ is the only vertex whose deletion from
$S^{\prime}$ produces a TDS, namely $S$. However, by symmetry, $C_{n}$ has
four $\gamma_{t}$-sets.

$(ii)\hspace{0.1in}$If $n\geq12$, then by Observation \ref{Ob_td_Paths} and
Proposition \ref{Prop_Gt_Cn}, $\Gamma_{t}(C_{n})\geq\gamma_{t}(C_{n})+2$.
Hence $C_{n}$ has a $P_{2}$-MTDS $S^{\ast}$ such that $|S^{\ast}|=2k+2$ and
$C_{n}[S^{\ast}]$ has $k+1$ components. Elementary calculations show that four
components are $P_{2}\overline{P}_{1}$ components.

$(iii)\hspace{0.1in}$By adding two vertices in succession, then deleting a
(different) vertex, we obtain%
\[
(\cdots\circ\circ\underline{\bullet}\bullet\overset{2}{\circ}%
\overset{2}{\bullet}\overset{2}{\circ}\overset{2}{\bullet}\cdots
)\rightarrow(\cdots\overset{2}{\circ}\underline{\bullet}%
\overset{2}{\blacklozenge}\overset{1}{\lozenge}\overset{3}{\blacklozenge
}\overset{1}{\lozenge}\overset{2}{\bullet}\cdots)\rightarrow(\cdots
\overset{2}{\circ}\underline{\bullet}\overset{2}{\bullet}\overset{2}{\lozenge
}\overset{2}{\blacklozenge}\overset{1}{\circ}\overset{2}{\bullet}\cdots),
\]
where \underline{$\bullet$} is a marker to indicate a specific vertex. Adding
and deleting another vertex, we obtain%
\[
(\cdots\overset{2}{\circ}\underline{\bullet}\overset{2}{\bullet}%
\overset{2}{\circ}\overset{2}{\bullet}\overset{1}{\circ}\overset{2}{\bullet
}\cdots)\rightarrow(\cdots\overset{2}{\circ}\underline{\bullet}%
\overset{2}{\bullet}\overset{2}{\circ}\overset{2}{\bullet}\overset{1}{\circ
}\overset{3}{\blacklozenge}\cdots)\rightarrow(\cdots\overset{2}{\circ
}\underline{\bullet}\overset{2}{\bullet}\overset{2}{\circ}\overset{2}{\bullet
}\overset{2}{\lozenge}\overset{2}{\blacklozenge}\cdots).
\]
Continuing the process, we eventually obtain the TDS $(\cdots
\overset{2}{\blacklozenge}\overset{1}{\lozenge}\underline{\bullet
}\overset{2}{\bullet}\overset{2}{\circ}\overset{2}{\bullet}\overset{2}{\circ
}\overset{2}{\bullet}\cdots)$ of cardinality $\gamma_{t}(C_{n})+1$, and one
last step -- a vertex deletion -- produces the $\gamma_{t}$-set
\[
S^{\prime}=(\cdots\overset{2}{\bullet}\circ\underline{\lozenge}%
\overset{2}{\bullet}\overset{2}{\circ}\overset{2}{\bullet}\overset{2}{\circ
}\overset{2}{\bullet}\cdots).
\]
Hence $S^{\prime}$ is obtained from $S$ by a rotation $v_{i}\rightarrow
v_{i+1}$ for each $i$. Repeating the procedure twice more shows that $S$ is
connected to each of the three other $\gamma_{t}$-sets of $C_{n}$.

Now let $S^{\ast}$ be a $P_{2}$-MTDS of cardinality $2k+2$. Then $n\geq12$ and
as shown in Lemma \ref{Lem_Cycle2}$(ii)$ we may assume that the four
$P_{2}\overline{P}_{1}$ components of $C_{n}[S^{\ast}]$ occur consecutively.
Hence $S^{\ast}$ is of the form $(\overset{2}{\bullet}\overset{1}{\circ
}\overset{2}{\bullet}\overset{1}{\circ}\overset{2}{\bullet}\overset{1}{\circ
}\overset{2}{\bullet}\overset{1}{\circ}\overset{2}{\bullet}\overset{2}{\circ
}\cdots)$. The operations%
\[
(\overset{2}{\bullet}\overset{1}{\circ}\overset{2}{\bullet}\overset{1}{\circ
}\overset{2}{\bullet}\overset{1}{\circ}\overset{2}{\bullet}\overset{1}{\circ
}\overset{2}{\bullet}\overset{2}{\circ}\cdots)\rightarrow(\overset{2}{\bullet
}\overset{1}{\circ}\overset{5}{\blacklozenge}\overset{1}{\circ}%
\overset{2}{\bullet}\overset{1}{\circ}\overset{2}{\bullet}\overset{2}{\circ
}\cdots)\rightarrow(\overset{2}{\bullet}\overset{2}{\lozenge}%
\overset{4}{\blacklozenge}\overset{1}{\circ}\overset{2}{\bullet}%
\overset{1}{\circ}\overset{2}{\bullet}\overset{2}{\circ}\cdots)\rightarrow
(\overset{2}{\bullet}\overset{2}{\circ}\overset{7}{\blacklozenge
}\overset{1}{\circ}\overset{2}{\bullet}\overset{2}{\circ}\cdots)\rightarrow
(\overset{2}{\bullet}\overset{2}{\circ}\overset{2}{\blacklozenge
}\overset{2}{\lozenge}\overset{2}{\blacklozenge}\overset{2}{\lozenge
}\overset{2}{\bullet}\overset{2}{\circ}\cdots)
\]
show that $S^{\ast}$ belongs to the same component of $D_{2k+3}^{t}(C_{n})$ as
a $\gamma_{t}$-set of $C_{n}$. The result follows by
transitivity.~$\blacksquare$

\bigskip

\noindent\textbf{Lemma \ref{Lem_Cycle4}}\emph{\hspace{0.1in}If }$3\leq n\leq
9$\emph{ and }$n\neq8$\emph{, then }$d_{0}(C_{n})=\Gamma_{t}(C_{n})+1$\emph{.}

\noindent\textbf{Proof.\hspace{0.1in}}The result is obvious for $n\in
\{3,4,5\}$ because $D_{3}^{t}(C_{3})\cong K_{1,3}$, $D_{3}^{t}(C_{4})\cong
C_{8}$ and $D_{4}^{t}(C_{5})\cong C_{10}$. All MTDS's of $C_{6}$ are of the
form $(\overset{2}{\bullet}\overset{1}{\circ}\overset{2}{\bullet
}\overset{1}{\circ})$ or $(\overset{4}{\bullet}\overset{2}{\circ})$, and one
easily obtains that $d_{0}(C_{6})=5=\Gamma_{t}(C_{6})+1$. All MTDS's of
$C_{7}$ are of the form $(\overset{2}{\bullet}\overset{1}{\circ}%
\overset{2}{\bullet}\overset{2}{\circ})$ and the result is easy to check.
Finally, all MTDS's of $C_{9}$ are of the form $(\overset{3}{\bullet
}\overset{2}{\circ}\overset{2}{\bullet}\overset{2}{\circ})$ and
$(\overset{2}{\bullet}\overset{1}{\circ}\overset{2}{\bullet}\overset{1}{\circ
}\overset{2}{\bullet}\overset{1}{\circ})$ and again the result follows
easily.~$\blacksquare$

\section{Problems for future work}

\label{Sec_Pr}

\begin{problem}
Determine $d_{0}(G)$ for other classes of graphs.
\end{problem}

\begin{problem}
Construct classes of graphs $G_{\alpha}$ such that the difference
$d_{0}(G_{\alpha})-\Gamma_{t}(G_{\alpha})\geq\alpha\geq2$ (or show that the
difference is bounded).
\end{problem}

\begin{problem}
Find more classes of graphs that can/cannot be realized as $k$-total
domination graphs.
\end{problem}

\begin{question}
When is $D_{k}^{t}(G)$ Hamiltonian?
\end{question}

\begin{question}
Which graphs $G$ satisfy $D_{k}^{t}(G)\cong G$ for some value of $k$?
\end{question}

\begin{question}
What is the complexity of determining whether two MTDS's of $G$ are in the
same component of $D_{k}^{t}(G)$, or of $D_{\Gamma(G)+1}^{t}(G)$?
\end{question}

\newpage

\section{Appendix:}

\section*{Upper total domination numbers of cycles}

Let $C_{n}=(v_{0},v_{1},...,v_{n-1},v_{0})$. When discussing subsets of
$\{v_{0},v_{1},...,v_{n-1}\}$ the arithmetic in the subscripts is performed
modulo $n$.

\begin{lemma}
\label{Lem_CnPn}For every integer $n\geq3$, $\Gamma_{t}(C_{n})\leq\Gamma
_{t}(P_{n})$.
\end{lemma}

\noindent\textbf{Proof.\hspace{0.1in}}Let $S$ be an MTDS of $C_{n}$. Then
$|S|<n$. Without loss of generality say $v_{0}\notin S$ and $v_{1}\in S$, and
let $P_{n}=C_{n}-v_{0}v_{n-1}=(v_{0},v_{1},...,v_{n-1})$. We show that $S$ is
an MTDS of $P_{n}$.

Note that $v_{1}$ totally dominates $v_{0}$ in $P_{n}$, and $S-\{v_{1}\}$ does
not dominate $v_{0}$ in $P_{n}$.

If $v_{n-1}\notin S$, then $v_{n-2}$ totally dominates $v_{n-1}$ in $C_{n}$ as
well as in $P_{n}$, and $S-\{v_{n-2}\}$ does not dominate $v_{n-1}$ in $P_{n}%
$. On the other hand, if $v_{n-1}\in S$, then $v_{n-2}\in S$ to totally
dominate $v_{n-1}$, and $v_{n-3}\notin S$, otherwise $S-\{v_{n-1}\}$ is a TDS
of $C_{n}$, contrary to the minimality of $S$. Hence $v_{n-2}$ is an isolated
vertex in the subgraph induced by $S-\{v_{n-1}\}$, which is therefore not a
TDS of $P_{n}$.

For any vertex $v\in S-\{v_{1},v_{n-2},v_{n-1}\}$, if $S-\{v\}$ totally
dominates $P_{n}$, then $S-\{v\}$ totally dominates $C_{n}$, which is not the
case. Therefore, $S$ is an MTDS of $P_{n}$. It follows that $\Gamma_{t}%
(C_{n})\leq\Gamma_{t}(P_{n})$.~$\blacksquare$

\bigskip

We next mention some obvious properties of minimal total dominating sets of
$C_{n}$.

\begin{remark}
Let $S$ be an MTDS of $C_{n}$. Then

\begin{enumerate}
\item[$(i)$] each component of $C_{n}[S]$ is either $P_{2},\ P_{3}$ or $P_{4}$;

\item[$(ii)$] each $P_{3}$ or $P_{4}$ component is preceded and followed by
exactly two consecutive vertices of $C_{n}-S$.
\end{enumerate}
\end{remark}

The following lemma will be used to describe the minimal total dominating sets
of $C_{n}$.

\begin{lemma}
\label{Lem_TDS}Let $S$ be an MTDS of $C_{n}$.

\begin{enumerate}
\item[$(i)$] If $C_{n}[S]$ has $P_{4}$ and $P_{2}$ components $(\cdots
\overset{2}{\circ}\overset{4}{\bullet}\overset{2}{\circ}\overset{2}{\bullet
}\circ\cdots)$, then $(\cdots\overset{2}{\circ}\overset{4}{\bullet
}\overset{2}{\circ}\overset{2}{\bullet}\circ\cdots)\rightarrow(\cdots
\overset{2}{\circ}\overset{2}{\bullet}\overset{1}{\circ}\overset{2}{\bullet
}\overset{1}{\circ}\overset{2}{\bullet}\circ\cdots)$ produces an MTDS
$S^{\prime}$ of $C_{n}$ such that $S^{\prime}$ has more $P_{2}$ components
than $S$, and $|S^{\prime}|=|S|$.

\item[$(ii)$] If $C_{n}[S]$ has only $P_{2}$ and $P_{4}$ components, then
$C_{n}$ has an MTDS $S^{\prime}$ such that $|S|=|S^{\prime}|$ and
$C_{n}[S^{\prime}]$ has only $P_{2}$ components.

\item[$(iii)$] If $C_{n}[S]$ has two consecutive $P_{3}$ components
$(\cdots\overset{2}{\circ}\overset{3}{\bullet}\overset{2}{\circ}%
\overset{3}{\bullet}\overset{2}{\circ}\cdots)$, then $(\cdots\overset{2}{\circ
}\overset{3}{\bullet}\overset{2}{\circ}\overset{3}{\bullet}\overset{2}{\circ
}\cdots)\rightarrow(\cdots\overset{2}{\circ}\overset{2}{\bullet}%
\overset{1}{\circ}\overset{2}{\bullet}\overset{1}{\circ}\overset{2}{\bullet
}\overset{2}{\circ}\cdots)$ produces an MTDS $S^{\prime}$ of $C_{n}$ with more
$P_{2}$ components than $S$, and $|S^{\prime}|=|S|$.

\item[$(iv)$] If $C_{n}[S]$ has only $P_{2}$ components, then $|S|\leq
2\left\lfloor \frac{n}{3}\right\rfloor $.
\end{enumerate}
\end{lemma}

\noindent\textbf{Proof.\hspace{0.1in}}Statements $(i),$ $(iii)$ and $(iv)$ are
obvious, while $(ii)$ follows from repeated applications of $(i)$ if
$C_{n}[S]$ has at least one $P_{2}$ component. Otherwise,
$S=(\overset{4}{\bullet}\overset{2}{\circ}\cdots\overset{4}{\bullet
}\overset{2}{\circ})$, and $S^{\prime}=(\overset{2}{\bullet}\overset{1}{\circ
}\overset{2}{\bullet}\overset{1}{\circ}\cdots\overset{2}{\bullet
}\overset{1}{\circ})$ is the desired set.~$\blacksquare$

\begin{proposition}
For any $n\geq3$,
\[
\Gamma_{t}(C_{n})=\left\{
\begin{tabular}
[c]{ll}%
$\left\lfloor \frac{2n}{3}\right\rfloor =\Gamma_{t}(P_{n})$ & if $n\equiv0$ or
$1\ (\operatorname{mod}\ 3)$\\
& \\
$\left\lfloor \frac{2n}{3}\right\rfloor <\Gamma_{t}(P_{n})$ & if
$n\equiv5\ (\operatorname{mod}\ 6)$\\
& \\
$2\left\lfloor \frac{n}{3}\right\rfloor <\left\lfloor \frac{2n}{3}%
\right\rfloor <\Gamma_{t}(P_{n})$ & if $n\equiv2\ (\operatorname{mod}\ 6).$%
\end{tabular}
\right.
\]

\end{proposition}

\noindent\textbf{Proof.\hspace{0.1in}}Let $C_{n}=(v_{0},v_{1},...,v_{n-1}%
,v_{0})$. If $n\equiv0\ (\operatorname{mod}\ 3)$, let $S=(\overset{2}{\bullet
}\overset{1}{\circ}\overset{2}{\bullet}\overset{1}{\circ}\cdots
\overset{2}{\bullet}\overset{1}{\circ})$. If $n\equiv1\ (\operatorname{mod}%
\ 3)$, let $S=(\overset{2}{\bullet}\overset{1}{\circ}\overset{2}{\bullet
}\overset{1}{\circ}\cdots\overset{2}{\bullet}\overset{1}{\circ}%
\overset{2}{\bullet}\overset{2}{\circ})$. In either case $S$ dominates $C_{n}$
and $|S|=2\left\lfloor \frac{n}{3}\right\rfloor =\left\lfloor \frac{2n}%
{3}\right\rfloor =2\left\lfloor \frac{n+1}{3}\right\rfloor =\Gamma_{t}(P_{n}%
)$. Moreover, $C_{n}[S]$ consists of $\left\lfloor \frac{n}{3}\right\rfloor $
disjoint copies of $K_{2}$, hence is an MTDS of $C_{n}$. Therefore $\Gamma
_{t}(C_{n})\geq|S|=\Gamma_{t}(P_{n})$ and the result follows from Lemma
\ref{Lem_CnPn}.

Assume next that $n\equiv5\ (\operatorname{mod}\ 6)$. The result is obvious
for $n=5$, hence say $n=6k+5$, $k\geq1$, and let $S=(\overset{4}{\bullet
}\overset{2}{\circ}\cdots\overset{4}{\bullet}\overset{2}{\circ}%
\overset{3}{\bullet}\overset{2}{\circ})$. Note that $S$ is an MTDS of $C_{n}$
and $|S|=4k+3=\left\lfloor \frac{2n}{3}\right\rfloor $. Hence $\Gamma
_{t}(C_{n})\geq\left\lfloor \frac{2n}{3}\right\rfloor $.

Let $X$ be a $\Gamma_{t}$-set of $C_{n}$ with as many $P_{2}$ components as
possible. Suppose first that $C_{n}[X]$ has no $P_{3}$ components. By Lemma
\ref{Lem_TDS}$(ii)$ we may assume that $C_{n}[X]$ has only $P_{2}$ components.
But then, by Lemma \ref{Lem_TDS}$(iv)$, $|X|\leq2\left\lfloor \frac{n}%
{3}\right\rfloor <\left\lfloor \frac{2n}{3}\right\rfloor =|S|$, which
contradicts $X$ being a $\Gamma_{t}$-set. Hence assume that $C_{n}[X]$ has a
$P_{3}$ component. By Lemma \ref{Lem_TDS}$(iii)$ we may assume that $C_{n}[X]$
has no consecutive $P_{3}$ components. If $C_{n}[X]$ has a $P_{3}$ component
between two $P_{2}$ components, $(\cdots\circ\overset{2}{\bullet
}\overset{2}{\circ}\overset{3}{\bullet}\overset{2}{\circ}\overset{2}{\bullet
}\circ\cdots)$, then $(\cdots\circ\overset{2}{\bullet}\overset{2}{\circ
}\overset{3}{\bullet}\overset{2}{\circ}\overset{2}{\bullet}\circ
\cdots)\rightarrow(\cdots\circ\overset{2}{\bullet}\overset{1}{\circ
}\overset{2}{\bullet}\overset{1}{\circ}\overset{2}{\bullet}\overset{1}{\circ
}\overset{2}{\bullet}\circ\cdots)$ produces a larger MTDS of $C_{n}$, which is
impossible. Hence each $P_{3}$ component is preceded or followed by a $P_{4}$ component.

Suppose $C_{n}[X]$ contains a segment $(\cdots\circ\overset{2}{\bullet
}\overset{2}{\circ}\overset{3}{\bullet}\overset{2}{\circ}\overset{4}{\bullet
}\overset{2}{\circ}\cdots)$. Then $(\cdots\circ\overset{2}{\bullet
}\overset{2}{\circ}\overset{3}{\bullet}\overset{2}{\circ}\overset{4}{\bullet
}\overset{2}{\circ}\cdots)\rightarrow(\cdots\circ\overset{2}{\bullet
}\overset{1}{\circ}\overset{2}{\bullet}\overset{1}{\circ}\overset{2}{\bullet
}\overset{2}{\circ}\overset{3}{\bullet}\overset{2}{\circ}\cdots)$ produces a
$\Gamma_{t}$-set with more $P_{2}$ components than $X$, contrary to the choice
of $X$. We have now proved that each $P_{3}$ component is preceded and
followed by $P_{4}$ components. If $n=11$, we have therefore shown that
$C_{n}[X]$ consists of one $P_{3}$ and one $P_{4}$ component, so we assume
that $n\geq17$.

If $C_{n}[X]$ contains a segment $(\cdots\overset{2}{\circ}\overset{3}{\bullet
}\overset{2}{\circ}\overset{4}{\bullet}\overset{2}{\circ}\overset{3}{\bullet
}\overset{2}{\circ}\cdots)$, then $(\cdots\overset{2}{\circ}%
\overset{3}{\bullet}\overset{2}{\circ}\overset{4}{\bullet}\overset{2}{\circ
}\overset{3}{\bullet}\overset{2}{\circ}\cdots)\rightarrow(\cdots
\overset{2}{\circ}\overset{2}{\bullet}\overset{1}{\circ}\overset{2}{\bullet
}\overset{1}{\circ}\overset{2}{\bullet}\overset{1}{\circ}\overset{2}{\bullet
}\overset{1}{\circ}\overset{2}{\bullet}\overset{2}{\circ}\cdots)$ produces a
$\Gamma_{t}$-set with more $P_{2}$ components than $X$, contrary to the choice
of $X$. A similar argument shows that any sequence of consecutive $P_{4}$
components of $C_{n}[X]$ cannot be preceded and followed by distinct $P_{3}$
components. Finally, by Lemma \ref{Lem_TDS}$(i)$ no $P_{4}$ component is
preceded or followed by a $P_{2}$ component. It follows that $C_{n}[X]$
consists of a single $P_{3}$ component, all other components being $P_{4}$'s.
By the choice of $X$ as a $\Gamma_{t}$-set with the largest number of $P_{2}$
components, we deduce that all $\Gamma_{t}$-sets of $C_{n}$ have this
property. Hence $\Gamma_{t}(C_{n})=4k+3=\left\lfloor \frac{2n}{3}\right\rfloor
$.

The final case to consider is $n\equiv2\ (\operatorname{mod}\ 6)$. Say
$n=6k+2$, $k\geq1$, and define $S=(\overset{2}{\bullet}\overset{1}{\circ
}\overset{2}{\bullet}\overset{1}{\circ}\overset{2}{\bullet}\cdots
\overset{1}{\circ}\overset{2}{\bullet}\overset{2}{\circ}\overset{2}{\bullet
}\overset{2}{\circ})$. Then $S$ is an MTDS of $C_{n}$ and $|S|=2\left\lfloor
\frac{6k-5+2}{3}\right\rfloor =4k-2=2\left\lfloor \frac{n}{3}\right\rfloor $,
hence $\Gamma_{t}(C_{n})\geq2\left\lfloor \frac{n}{3}\right\rfloor $.

Suppose $X$ is a $\Gamma_{t}$-set of $C_{n}$ with as many $P_{2}$ components
as possible that also contains a $P_{3}$ component. Arguing as above we come
to the same conclusion, namely that $C_{n}[X]$ consists of a single $P_{3}$
component, all other components being $P_{4}$'s. But then a simple counting
argument shows that $n=6t+5$ for some integer $t$, which is not the case since
$n\equiv2\ (\operatorname{mod}\ 6)$. Hence any $\Gamma_{t}$-set of $C_{n}$
contains only $P_{2}$ and $P_{4}$ components. By Lemma \ref{Lem_TDS}$(ii)$ we
may thus assume that $C_{n}$ has a $\Gamma_{t}$-set that has only $P_{2}$
components, so that, by Lemma \ref{Lem_TDS}$(iv)$, $\Gamma_{t}(C_{n}%
)\leq2\left\lfloor \frac{n}{3}\right\rfloor $ and we are done.~$\blacksquare$


\begin{thebibliography}{99}                                                                                               %


\bibitem {davood}S.~Alikhani, D.~Fatehi and S.~Klav\v{z}ar, \textit{On the
structure of dominating graphs}, Graphs Combin.\ \textbf{33} (2017), 665--672.

\bibitem {BC}P.~Bonsma, L.~Cereceda, Finding paths between graph colourings:
PSPACE-completeness and superpolynomial distances, \emph{Theor. Comput. Sci.}
\textbf{410 }(50) (2009), 5215--5226.

\bibitem {CHJ1}L.~Cereceda, J.~van den Heuvel, M.~Johnson, Connectedness of
the graph of vertex-colourings, \emph{Discrete Math.} \textbf{308} (2008), 913--919.

\bibitem {CHJ2}L.~Cereceda, J.~van den Heuvel, M.~Johnson, Finding paths
between 3-colorings, \emph{J. Graph Theory} \textbf{67} (2011), 69--82.

\bibitem {CDH}E.~J.~Cockayne, R.~M.~Dawes, S.~T.~Hedetniemi, Total domination
in graphs, \emph{Networks} \textbf{10} (1980), 211--219.

\bibitem {Dorbec}P.~Dorbec, M.~A.~Henning, J.~McCoy, Upper total domination
versus upper paired domination, \emph{Quaest. Math}. \textbf{30} (2007), 1--12.

\bibitem {FH}O.~Favaron, M.\ A.\ Henning, Upper total domination in claw-free
graphs. \emph{J. Graph Theory} \textbf{44} (2003), 148--158.

\bibitem {HS1}R.~Haas, K.~Seyffarth, The $k$-dominating graph, \emph{Graphs
Combin.} \textbf{30 }(3) (2014), 609--617.

\bibitem {HS2}R.~Haas, K.~Seyffarth, Reconfiguring dominating sets in some
well-covered and other classes of graphs, \emph{Discrete Math.} \textbf{340
}(2017), 1802--1817.

\bibitem {haddadan-2015}A.~Haddadan, T.~Ito, A.~E.~Mouawad, N.~Nishimura,
H.~Ono, A.~Suzuki, Y.~Tebbal, \textit{The complexity of dominating set
reconfiguration}, \emph{Theoret. Comput. Sci.} \textbf{651} (2016), 37--49.

\bibitem {HHS}T.~W.~Haynes, S.~T.~Hedetniemi, P.~J.~Slater, \emph{Fundamentals
of Domination in Graphs}. Marcel Dekker, New York, 1998.

\bibitem {HY}M.~A.~Henning, A.~Yeo, \emph{Total Domination in Graphs}.
Springer, New York, 2014.

\bibitem {Ito2}T.~Ito, E.~D.~Demaine, N.~J.~A.~Harvey, C.~H.~Papadimitriou,
M.~Sideri, R.~Uehara, Y.~Uno, On the complexity of reconfiguration problems,
\emph{Theoret. Comput. Sci.} \textbf{412}(12--14) (2011), 1054--1065.

\bibitem {IKD}T.~Ito, M.~Kaminski, E.~D.~Demaine, Reconfiguration of list
edge-colorings in a graph, in: \emph{Algorithms and data structures, LNCS},
vol. 5664, Springer, Berlin, 2009, 375--386.

\bibitem {MRT}C.~M.~Mynhardt, A.~Roux,\ L.~E.~Teshima, Connected
$k$-dominating graphs, arXiv:1708.05458v2 [math.CO].

\bibitem {SMN}A.~Suzuki, A.~E.~Mouawad, N.~Nishimura, Reconfiguration of
dominating sets, in: Z.~Cai, A.~Zelikovsky, A.~Bourgeois (Eds.), COCOON 2014,
\emph{LNCS}, vol. 8591, Springer, Heidelberg, 2014, 405--416.
\end{thebibliography}
\end{document}